\documentclass[journal]{IEEEtran}
\ifCLASSINFOpdf
\usepackage[pdftex]{graphicx}
\else
\fi
\usepackage{amsmath,amssymb}
\usepackage {autobreak}
\newcommand{\CSizeVec}[1]{ \in {\mathbb C}^{#1} }
\newcommand{\CSize}[2]{ \in {\mathbb C}^{#1 \times #2} }

\newcommand{\Field}[1]{ \mbox{\boldmath $#1$} }
\def\figuref{Fig.~}

\def\tableref{TABLE~}
\usepackage[hang,small,bf]{caption}
\usepackage[subrefformat=parens]{subcaption}
\begin{document}
%
\title{Characteristic Basis Function Method Combined with Calder\'{o}n Multiplicative Preconditioner for PMCHWT Formulation}
%
%
%

\author{Tai~Tanaka,~\IEEEmembership{Member,~IEEE,}
        Kazuki~Niino, ~\IEEEmembership{Member,~IEEE,}
        and~Naoshi~Nishimura\IEEEmembership{}
\thanks{T. Tanaka, K. Niino and N. Nishimura are with the Kyoto University, Kyoto, Japan.}
}

\maketitle

\begin{abstract}
We propose a novel characteristic basis function method for analyzing the scattering by dielectric objects based on the Poggio-Miller-Chang-Harrington-Wu-Tsai formulation. 
In the proposed method, the electric and magnetic currents are orthogonalized with the help of the singular value decomposition, and are used as dual basis functions in a way similar to the RWG and BC basis functions. We show that the use of the Calder\'{o}n multiplicative precondtioner together with the proposed method can prevent from the poor convergence of the solution of the matrix equation in problems involving dielectrics.
We considered three different shapes of dielectric scatterers for the purpose of validation. 
The numerical results agreed well with those obtained by the conventional method of moments and the proposed method was faster than the conventional method. 
These results indicate that the proposed method is effective for scattering analysis of the dielectrics. 
\end{abstract}

\begin{IEEEkeywords}
Method of Moments~(MoM), Characteristic Basis Function Method~(CBFM), Poggio-Miller-Chang-Harrington-Wu-Tsai~(PMCHWT) formulation, Calder\'{o}n Multiplicative Preconditioner~(CMP). 
\end{IEEEkeywords} 

%
\IEEEpeerreviewmaketitle
\section{Introduction}
\IEEEPARstart{T}{he} method of moments~(MoM) \cite{Harrington1993} is one of numerical methods of analyzing electromagnetic field problems which are based on boundary integral equations. 
In the MoM one utilizes the Galerkin method \cite{Harrington1993} to discretize integral equations into a system of linear equations, which is solved with direct or iterative algorithms \cite{Saad2003}. 
The computational cost of the MoM is in general governed by that of computing the impedance matrix in the discretized system, which is usually quite expensive since the impedance matrix is dense. 
One of solutions to this problem is the use of fast methods such as fast multipole method~(FMM), which accelerate the calculation of products of the impedance matrix with given vectors.
Hence iterative algorithms accelerated with the FMM are widely used for analyzing large-scale problems \cite{Chew2001}--\cite{Song1997}.
This paper focuses on such iterative algorithms, although fast methods for direct solvers such as those based on the $\mathcal{H}$-matrix have been developed recently \cite{Grasedyck2003}--\cite{Gholami2018}.

Analyses using iterative methods often suffer from poor convergence, which is usually solved by preconditioners.
In recent years, the Calder\'{o}n multiplicative preconditioner~(CMP) has been proposed in order to improve the convergence of the electric field integral equation~(EFIE) \cite{Andriulli2008}.
This preconditioner is based on simple mathematical relations between products of integral operators.
However, discretizing products of integral operators in Maxwell's equations is more complicated than is expected from the apparent mathematical simplicity.
In the case of the electric field integral equation for example, two types of mutually (almost) orthogonal basis functions for the same functional space are necessary for discretizing operators corresponding to the preconditioner and the impedance matrix.
Indeed, it is known that a naive choice of standard basis functions such as the Rao-Wilton-Glisson~(RWG) basis function \cite{Rao1982}, denoted by $\Field{f}$, for discretizing both the impedance and preconditioning matrices does not work since the Gram matrix associated with this choice is singular, while its inverse appears in the formulation of the CMP.
Using RWG for the impedance matrix and its $90$ degree rotation, i.e.\ $\hat{\Field{n}}\times{\Field{f}}$, for the preconditioner does not make sense since $\hat{\Field{n}}\times\Field{f}$ does not belong to $H_{\mathrm{div}}$ where $\hat{\Field{n}}$ is the unit normal vector on the boundary.
As a matter of fact, 
a widely used choice of basis functions for the CMP is the RWG basis function for the impedance matrix and the Buffa-Christiansen~(BC) basis function \cite{Buffa2007} for the preconditioner.
However the computational time of the preconditioner discretized with the BC function is much more than that of the impedance matrix with the RWG function since the BC function is defined on the barycentric refinement of an original mesh on which the RWG function is defined.
The CMP has also been applied to integral equations for scattering problems with dielectrics such as the Poggio-Miller-Chang-Harrington-Wu-Tsai (PMCHWT) formulation \cite{Yan2010, Cools2011}.
For the PMCHWT, another formulation of the CMP has been proposed \cite{Niino2012, yla2012stable}, in which the surface electric and magnetic currents are expanded with the RWG and BC basis functions, respectively.
This formulation has a beautiful symmetry and is able to accelerate the convergence in problems for scatterers having smooth boundaries,
but the accuracy of the surface current spanned with the BC basis function is usually worse than that with the RWG function when the boundary contains sharp edges or corners \cite{Niino2012}.

As another class of numerical methods to speed up MoM apart from the preconditioned fast methods, one may mention domain decomposition methods \cite{Heldrig2007, Matekovits2007}, among which we are particularly interested in the characteristic basis function method~(CBFM) \cite{Prakash2003}--\cite{Chen2015}.
The CBFM is known to be particularly suited for scattering problems by finite periodic scatterers such as radiation by array antennas, and is also known to be quite effective in calculating monostatic radar cross sections (RCS).
The CBFM first decomposes domains into several pieces called cells, and roughly solves scattering problems (which may be called generating problems) in each cell with certain numbers of incident waves. 
After obtaining as many solutions as the number of the incident waves in each cell, we generate
linearly independent basis functions from these solutions by applying the singular value decomposition~(SVD) to a matrix in which these solutions are arranged as column vectors.
The basis functions obtained in this way are called the characteristic basis functions~(CBFs) and the CBFM solves the original problem by using CBFs as the basis functions.
The computational time of the CBFM is less than that of the ordinary MoM since the number of CBFs is much less than that of the standard basis functions.
The accuracy of the CBFM depends significantly on the choice of CBFs.
We have proposed a class of CBFs, called improved primary CBFs~(IPCBFs), which are expected to approximate surface currents efficiently \cite{Tanaka2016}--\cite{Tanaka2019} because they are more accurate solutions of the generating problems than the standard (primary) CBFs. 
The IPCBFs take into account the influence of higher-order CBFs to primary CBFs iteratively.
The accuracy of IPCBFs can be controlled as one sets the residual norm of the iterative CBF generations appropriately.
In CBFM one often uses direct methods for solving linear systems since the CBFM considerably reduces the number of unknowns.
In large-scale problems, however, the number of unknowns is not always sufficiently small for using direct solvers.
Also the convergence of iterative algorithms in the CBFM can be as poor as that in the standard MoM since they are based on the same ill-conditioned integral equations.
It is therefore very important to reduce the number of iterations in the iterative CBFMs.

The goal of this paper is to propose a CMP in the CBFM.
We numerically construct two types of mutually-orthogonal basis functions required in the CMP
by applying the SVD to the Gram matrix of basis functions obtained in the process of the CBFM.
We will show that the basis functions obtained in this way enable the use of the symmetric CMP in the PMCHWT formulation,
in which these two types of basis functions, respectively, expand the surface electric and magnetic currents, without suffering from bad accuracy due to sharp edges or corners on the boundary observed in \cite{Niino2012}.
Also this formulation seems to be more natural from the aspect of the CBFM since CBFs in this method are bi-orthogonal as vector fields
while the conventional CBFs are orthogonalized as algebraic vectors.

The rest of this paper consists of the following four parts. 
In the next section, we give a preliminary description of the integral equations and the MoM used in this paper.  
Section \ref{sec_cbfm_with_cmp} describes the new CBFM with the CMP.
In section \ref{sec_num_ex} we validate the formulation by analyzing three types of scatterers. 
Finally we give some concluding remarks and discuss future works in section \ref{sec_conclusion}.

In this paper we use the following definitions and symbols. 
RCSs will refer to monostatic RCSs according to the definition in \cite{Knott1993}. 
Italicized bold symbols such as $\Field{A}$ represent three-dimensional vector fields. 
Capital and small symbols having bold and upright type, such as ${\bf A}$ and ${\bf a}$, represent matrices and vectors, respectively. 
Matrix and vector symbols enclosed in square brackets with subscripts represent elements of matrices or vectors, such as $\begin{bmatrix} {\bf A} \end{bmatrix}_{ij}$. 
A matrix or vector enclosed in round brackets with subscripts, such as $\left({\bf A}\right)_{ij}$, represents a submatrix or a subvector, respectively. 
The inner product $\Field{A} \cdot \Field{B}$ is defined as $\sum_i \bar{a}_i b_i$, where $a_i$, $b_i$ are the components of vectors $\Field{A}$ and $\Field{B}$, and $\bar{a}_i$ is the complex conjugate of $a_i$. 

\section{Preliminary Formulations}
In this section, we formulate the electromagnetic scattering problem for homogeneous dielectric objects and MoM as a method for solving it. 
\subsection{PMCHWT Formulation}
\begin{figure}[!t]
  \centering
      \includegraphics[keepaspectratio, scale=0.4]{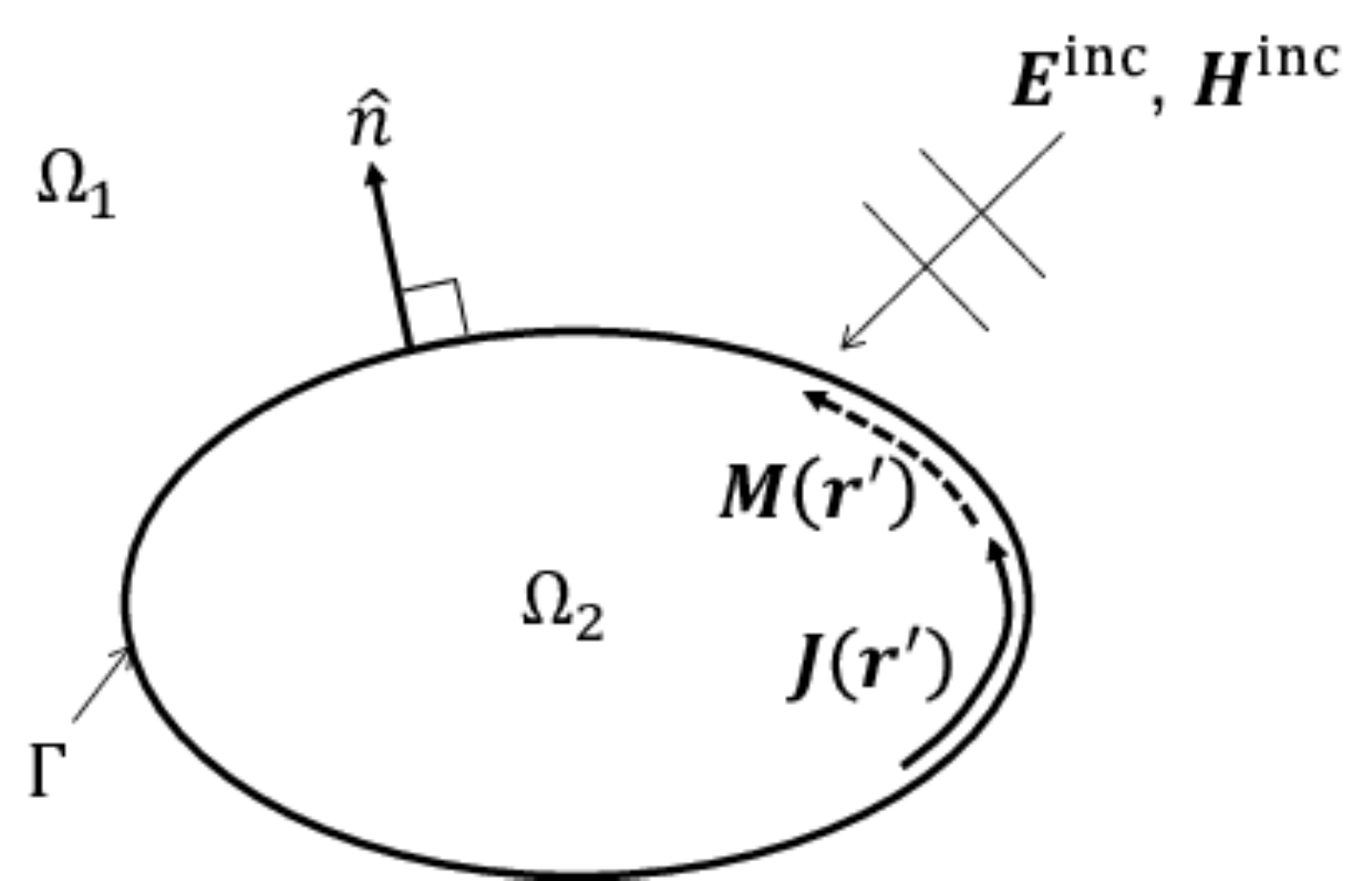}
  \caption{Definition of the domain}
  \label{fig_domain}
\end{figure}

We consider an electromagnetic wave scattering problem with a homogeneous dielectric object shown in \figuref\ref{fig_domain}.
For simplicity we assume that we have a single scatterer in this section.
Extension to multiple scatterers is straightforward and indeed we will show some numerical examples with multiple scatterers in section \ref{sec_num_ex}.
The domains outside and inside of the scatterer are denoted by $\Omega_1$ and $\Omega_2=\mathbb{R}^{3} \setminus \overline{\Omega_{1}}$, respectively. 
The boundary $\Gamma = \partial \Omega_{2}$ is a closed surface and the unit normal vector $\hat{\Field{n}}$ points outward.
We assume the time factor to be $e^{j\omega t}$ with the frequency $\omega$.

The incident electromagnetic fields at a point $\Field{r}$ are $\Field{E}^{\rm inc}\left(\Field{r}\right)$ and $\Field{H}^{\rm inc}\left(\Field{r}\right)$, respectively.
The induced electromagnetic currents on the dielectric boundary $\Gamma$ are defined as $\Field{J}\left(\Field{r}\right)$ and $\Field{M}\left(\Field{r}\right)$.
For this wave scattering problem the integral equations based on the PMCHWT formulation \cite{Chew2001, Ergul2014} are written as follows:
\begin{align}
  \label{eq_PMCHWT1}
  \sum_{i=1}^{2}\left\{\left(\mathcal{K}_{i} \Field{M} \right) \left(\Field{r}\right) - \eta_{i} \left(\mathcal{T}_{i}\Field{J} \right)\left(\Field{r}\right) \right\} &= \hat{\Field{n}} \times \Field{E}^{\rm inc}\left(\Field{r}\right) \\
  \label{eq_PMCHWT2}
  \sum_{i=1}^{2}\left\{\frac{i}{\eta_{i}}\left(\mathcal{T}_{1} \Field{M} \right)\left(\Field{r}\right) + \left( \mathcal{K}_{i}\Field{J} \right)\left(\Field{r}\right)\right\}  &= - \hat{\Field{n}} \times \Field{H}^{\rm inc}\left(\Field{r}\right),
\end{align}
where $\epsilon_i$, $\mu_i$, and $\eta_{i}=\sqrt{\mu_{i}/\epsilon_{i}}$ are the permittivity, permeability, and wave impedance in $\Omega_i$, respectively.
The integral operators $\mathcal{K}_{i}$ and $\mathcal{T}_{i}$ are defined by
\begin{align}
  \label{eq_operator_K}
  {}\left(\mathcal{K}_{i} \Field{X} \right)\left(\Field{r}\right) = \hat{\Field{n}} \times P.V. \int_{\Gamma} \nabla G_{i}\left(\Field{r},\Field{r}^{\prime} \right) \times \Field{X} \left(\Field{r}^{\prime}\right) d \Field{r}^{\prime}& \\
  \label{eq_operator_T}
  {}\left(\mathcal{T}_{i} \Field{X} \right)\left(\Field{r}\right) = j k_{i} \hat{\Field{n}} \times F.P.\int_{\Gamma} \left( \mathcal{I} + \frac{\nabla\nabla}{k_{i}^{2}}\right)G_{i}\left(\Field{r},\Field{r}^{\prime}  \right) &\notag\\
\cdot \Field{X} \left(\Field{r}^{\prime}\right) d\Field{r}^{\prime}&
\end{align}
where $G_{i}\left(\Field{r},\Field{r}^{\prime} \right) $ is Green's function of the Helmholtz equation for the observation point $\Field{r}$ and the source point $\Field{r}^{\prime}$ in three dimension represented by
\begin{align}
  G_{i}\left(\Field{r},\Field{r}^{\prime}  \right) = \frac{e^{-j k_{i} |\Field{r}-\Field{r}^{\prime}|}}{4\pi |\Field{r}-\Field{r}^{\prime}|}.
\end{align}
In \eqref{eq_operator_K} and \eqref{eq_operator_T}, $P.V.$, $F.P.$ and $\mathcal{I}$ stand for the Cauchy principal value integral, finite part and the identity operator, respectively. 
The parameter $k_i$ is the wave number in $\Omega_i$ defined as $k_{i}=2\pi/\lambda_{i}$ with the wavelength $\lambda_i$.
Omitting the notation ``$\left(\Field{r}\right)$'' for simplicity, we can express the integral equations in \eqref{eq_PMCHWT1} and \eqref{eq_PMCHWT2} by
\begin{align}
  \label{eq_PMCHWT_matrix}
  \begin{bmatrix}
    \sum_{i=1}^{2}\mathcal{K}_{i}                           & -\sum_{i=1}^{2}\eta_{i} \mathcal{T}_{i} \\
    \sum_{i=1}^{2}\frac{1}{\eta_{i}}\mathcal{T}_{i}  & \sum_{i=1}^{2}\mathcal{K}_{i} \\
  \end{bmatrix}
  \begin{bmatrix}
    \Field{M} \\
    \Field{J}
  \end{bmatrix}
  =
  \begin{bmatrix}
    \hat{\Field{n}} \times \Field{E}^{\rm inc} \\
    - \hat{\Field{n}} \times \Field{H}^{\rm inc}
  \end{bmatrix}.
\end{align}

\subsection{Discretization of the Integral Equation}
The PMCHWT formulation in \eqref{eq_PMCHWT_matrix} is discretized with the Galerkin method \cite{Harrington1993} in this paper. 
The electromagnetic currents are approximated by linear combinations of $H_{\rm div}^{-\frac12}(\Gamma)$-conforming basis functions
$\Field{b}^{\rm J,M}_{n}\left(\Field{r}\right)$ with the expansion coefficients $\alpha_{n}^{\rm J,M}$ as follows:
\begin{align}
  \label{eq_J_expand}
  \Field{J}\left(\Field{r}\right) &\approx \sum_{n=1}^{N^{{\rm J}}}\alpha^{\rm J}_{n} \Field{b}_{n}^{\rm J}\left(\Field{r}\right) \\
  \label{eq_M_expand}
  \Field{M}\left(\Field{r}\right) &\approx \sum_{n=1}^{N^{{\rm M}}}\alpha^{\rm M}_{n} \Field{b}_{n}^{\rm M}\left(\Field{r}\right) 
\end{align}
where $N^{{\rm J}}$ and $N^{{\rm M}}$ are the total numbers of the basis functions for the electric and magnetic currents, denoted by $\Field{b}_{n}^{\rm J}$ and $\Field{b}_{n}^{\rm M}$, respectively.
Hence the total number of the unknowns $N$ for our problems is $N^{\rm J}+N^{\rm M}$.
Since the scatterer defined in this paper does not contain perfect electric conductors, $N^{\rm J}$ is equal to $N^{\rm M}$.
The discretized operators $\left[{\bf K}_{i}\right]_{m,n}$, $\left[{\bf T}_{i}\right]_{m,n}$ for elements $m$ and $n$ are described as 
\begin{align}\label{eq_K}
  \left[{\bf K}_{i}^{tb}\right]_{m,n} =& \int_{\Field{t}_{m}} \hat{\Field{n}} \times\Field{t}_{m}\left(\Field{r}\right) {\cdot}\nonumber \\  \left( \hat{\Field{n}} \times P.V.\vphantom{\int_{\Field{b}_{n}}}\right.& \left.\int_{\Field{b}_{n}} \nabla G_{i}\left( \Field{r},\Field{r}^{\prime}  \right) \times \Field{b}_{n} \left( \Field{r}^{\prime} \right) d \Field{r}^{\prime} \right) d \Field{r}, \\\label{eq_T}
\left[{\bf T}_{i}^{tb}\right]_{m,n} =& j k_{i} \int_{\Field{t}_{m}} \hat{\Field{n}} \times\Field{t}_{m}\left(\Field{r}\right){\cdot} \nonumber \\  \left( \hat{\Field{n}} \times F.P.\vphantom{\int_{\Field{b}_{n}}}\right.&\left.\int_{\Field{b}_{n}} \left( \mathcal{I} + \frac{\nabla\nabla}{k_{i}^{2}}\right)G_{i}\left(\Field{r},\Field{r}^{\prime}  \right) \cdot {\Field{b}_{n}} \left(\Field{r}^{\prime} \right) d \Field{r}^{\prime} \right) d \Field{r}
\end{align}
where $\Field{t}_{n}$ is a test function and $\int_{\Field{t}_{n}}$ stands for an integral over the support of $\Field{t}_{n}$. 
The superscripts $t$ and $b$ in ${\bf K}_{i}^{tb}$ and ${\bf T}_{i}^{tb}$ represent the types of test and trial (basis) functions. 
The discretized incident electric field $\left[{\bf v}^{{\rm E}}_{t}\right]_{m}$ and magnetic field $\left[{\bf v}^{{\rm H}}_{t}\right]_{m}$ corresponding to element $m$ are expressed by
\begin{align*}
  \left[{\bf v}^{{\rm E}}_{t}\right]_{m} &= \int_{\Field{t}_{m}} \hat{\Field{n}} \times\Field{t}_{m}\left(\Field{r}\right) \cdot \left( \hat{\Field{n}} \times \Field{E}^{\rm inc}\left(\Field{r}\right) \right) d \Field{r},\\
  \left[{\bf v}^{{\rm H}}_{t}\right]_{m} &= -\int_{\Field{t}_{m}} \hat{\Field{n}} \times\Field{t}_{m}\left(\Field{r}\right) \cdot \left( \hat{\Field{n}} \times \Field{H}^{\rm inc}\left(\Field{r}\right) \right) d \Field{r}.
\end{align*}
With these notations, the discretized matrix equation can be written as follows.
\begin{align}
  \label{eq_MoM}
  {\bf Z}_{t^Jt^Mb^Mb^J} {\bf j} &= {\bf v}_{t^Jt^M},\\
  \label{eq_Z}
  {\bf Z}_{t^Jt^Mb^Mb^J}
  &=\sum_{i=1}^{2}
  \begin{bmatrix}
    {\bf K}_{i}^{t^Jb^M}                          & -\eta_{i} {\bf T}_{i}^{t^Jb^J} \\
    \frac{1}{\eta_{i}}{\bf T}_{i}^{t^Mb^M} & {\bf K}_{i}^{t^Mb^J} \\
  \end{bmatrix} \CSize{N}{N},\\
  \label{eq_v}
  {\bf v}_{t^Jt^M}
  &=
  \begin{bmatrix}
    {\bf v}^{{\rm E}}_{t^J} \\
    {\bf v}^{{\rm H}}_{t^M}
  \end{bmatrix} \CSizeVec{N},
\end{align}
and the expansion coefficient vector ${\bf j}\CSizeVec{N}$ is written as 
\begin{align}
  \label{eq_coef_vector}
  {\bf j}
  =
  \begin{bmatrix}
    \begin{bmatrix}
      \alpha^{\rm M}_{1} & \cdots & \alpha^{\rm M}_{N^{{\rm M}}}
    \end{bmatrix} & 
    \begin{bmatrix}
      \alpha^{\rm J}_{1} & \cdots & \alpha^{\rm J}_{N^{{\rm J}}} 
    \end{bmatrix}
  \end{bmatrix}^{T},
\end{align}
where the superposed $T$ stands for the transpose of a vector. 
The 1st and 2nd subscripts of ${\bf Z}_{t^Jt^Mb^Mb^J}$ respectively correspond to the types of test functions for the 1st and 2nd row of \eqref{eq_Z} while the 3rd and 4th subscripts to the types of trial functions for the surface magnetic and electric currents.
Also the 1st and 2nd subscripts of ${\bf v}_{t^Jt^M}$ represent the types of test functions for the 1st and 2nd rows of \eqref{eq_v}, respectively. 
Note that we will eventually use different basis functions for the test functions of the 1st and 2nd rows in \eqref{eq_Z} and \eqref{eq_v} as well as for the trial functions of the magnetic and electric currents.
This is crucial in the use of the CMP as will be shown in the next section.

\subsection{Calder\'{o}n Preconditioning} \label{sec_pre-cald}
The Calder\'{o}n preconditioning is one of numerical methods to accelerate the convergence of iterative linear solvers for integral equations.
This method constructs a preconditioner based on the Calder\'{o}n formulae given by
\begin{align*}
  \begin{bmatrix}
    \mathcal{K}_i & -\eta_i \mathcal{T}_i \\
    \frac{1}{\eta_i} \mathcal{T}_i & \mathcal{K}_i \\
  \end{bmatrix}
  \begin{bmatrix}
    \mathcal{K}_i & -\eta_i \mathcal{T}_i \\
    \frac{1}{\eta_i} \mathcal{T}_i & \mathcal{K}_i \\
  \end{bmatrix}
  =
  \frac{\mathcal{I}}{4}.
\end{align*}
From this equation we expect that the square of the operator in \eqref{eq_PMCHWT_matrix} is well-conditioned.
Indeed it is true in a sense that the operator satisfies
\begin{align*}
  &\begin{bmatrix}
    \sum_{i=1}^{2}\mathcal{K}_{i}                    & - \sum_{i=1}^{2} \eta_{i} \mathcal{T}_{i} \\
    \sum_{i=1}^{2}\frac{1}{\eta_{i}}\mathcal{T}_{i}  & \sum_{i=1}^{2}\mathcal{K}_{i}  \\
   \end{bmatrix}\\
  &\cdot 
  \begin{bmatrix}
    \sum_{i=1}^{2}\mathcal{K}_{i}                    & - \sum_{i=1}^{2} \eta_{i} \mathcal{T}_{i} \\
    \sum_{i=1}^{2}\frac{1}{\eta_{i}}\mathcal{T}_{i}  & \sum_{i=1}^{2}\mathcal{K}_{i}  \\
  \end{bmatrix}
  = \mathcal{S} + \mathcal{K}
\end{align*}
where $\mathcal{S}$ is a bounded operator with a bounded inverse and $\mathcal{K}$ is a compact operator.
Hence we can construct a right preconditioner by discretizing the integral equation
\begin{align}
  \nonumber &\begin{bmatrix}
    \sum_{i=1}^{2}\mathcal{K}_{i}                    & - \sum_{i=1}^{2} \eta_{i} \mathcal{T}_{i} \\
    \sum_{i=1}^{2}\frac{1}{\eta_{i}}\mathcal{T}_{i}  & \sum_{i=1}^{2}\mathcal{K}_{i}  \\
  \end{bmatrix}\\ \label{eq_precond_pmchwt_matrix}
  &\cdot
  \begin{bmatrix}
    \sum_{i=1}^{2}\mathcal{K}_{i}                    & - \sum_{i=1}^{2} \eta_{i} \mathcal{T}_{i} \\
    \sum_{i=1}^{2}\frac{1}{\eta_{i}}\mathcal{T}_{i}  & \sum_{i=1}^{2}\mathcal{K}_{i}  \\
  \end{bmatrix}
  \begin{bmatrix}
    \Field{M}^{\prime} \\
    \Field{J}^{\prime}
  \end{bmatrix}
  =
  \begin{bmatrix}
    \hat{\Field{n}} \times \Field{E}^{\rm inc} \\
    - \hat{\Field{n}} \times \Field{H}^{\rm inc}
  \end{bmatrix},\\\nonumber
  \intertext{where}
  &
  \begin{bmatrix}
    \Field{M} \\
    \Field{J}
  \end{bmatrix}
  =
  \begin{bmatrix}
    \sum_{i=1}^{2}\mathcal{K}_{i}                    & - \sum_{i=1}^{2} \eta_{i} \mathcal{T}_{i} \\
    \sum_{i=1}^{2}\frac{1}{\eta_{i}}\mathcal{T}_{i}  & \sum_{i=1}^{2}\mathcal{K}_{i}  \\
  \end{bmatrix}
  \begin{bmatrix}
    \Field{M}^{\prime} \\
    \Field{J}^{\prime}
  \end{bmatrix}.
\end{align}
It is known that a naive discretization of equation \eqref{eq_precond_pmchwt_matrix} which uses only the RWG basis functions $\Field{f}$ for both test and trial functions (for both $\bf J$ and $\bf M$)
\begin{align*}
  &{\bf Z}_{ffff}{\bf G}_{ffff}^{-1} {\bf Z}_{ffff}{\bf G}_{ffff}^{-1} \tilde{\bf j}= {\bf v}_{ff},\\
  &\quad {\bf j}={\bf G}_{ffff}^{-1} {\bf Z}_{ffff}{\bf G}_{ffff}^{-1}\tilde{\bf j},
\end{align*}
does not make sense since the Gram matrix of the RWG basis function ${\bf j}$ defined by
\begin{align}
  \nonumber
  &{\bf G}_{ffff} \\
  \label{eq_Gffff}  &=\begin{bmatrix}
    \int_\Gamma(\hat{\Field{n}}\times \Field{f}_i({\bf r}))\cdot \Field{f}_j({\bf r}) d\Field{r} & 0 \\
    0 &  \int_\Gamma(\hat{\Field{n}}\times \Field{f}_i({\bf r}))\cdot \Field{f}_j({\bf r}) d\Field{r}
  \end{bmatrix}
\end{align}
is a singular matrix,
where ${\bf Z}_{ffff}$ is the matrix defined in \eqref{eq_Z} with the RWG basis functions for all the trial and test functions.
One of possible discretization methods to avoid this problem is to use the BC basis functions $\Field{g}$ for the preconditioner, namely discretizing \eqref{eq_precond_pmchwt_matrix} into
\begin{align*}
  &{\bf Z}_{ffff} {\bf G}_{ggff}^{-1} {\bf Z}_{gggg} {\bf G}_{ffgg}^{-1} \tilde{\bf j}= {\bf v}_{ff},\\
  &\quad {\bf j} = {\bf G}_{ggff}^{-1} {\bf Z}_{gggg} {\bf G}_{ffgg}^{-1} \tilde{\bf j}, 
\end{align*}
where
\begin{align*}
  &{\bf G}_{ggff} =  -{\bf G}_{ffgg} \\
  =&
  \begin{bmatrix}
    \int_\Gamma\left(\hat{\Field{n}}\times \Field{g}_i\left(\Field{r}\right)\right)\cdot \Field{f}_j\left(\Field{r}\right) d\Field{r} & 0 \\
    0 & \int_\Gamma\left(\hat{\Field{n}}\times \Field{g}_i\left(\Field{r}\right)\right)\cdot \Field{f}_j\left(\Field{r}\right) d\Field{r}
  \end{bmatrix},
\end{align*}
and ${\bf Z}_{gggg}$ is the matrix in \eqref{eq_Z} with the BC basis functions used for all the trial and test functions.
Another alternative is to expand the electric and magnetic currents with the RWG and BC basis functions, respectively, with appropriately chosen test functions \cite{yla2012stable, Niino2012}, i.e.:
\begin{align}
  \label{eq_precond_pmchwt_mixed_matrix}
  &{\bf Z}_{gffg} {\bf G}_{gffg}^{-1} {\bf Z}_{gffg} {\bf G}_{gffg}^{-1} \tilde{\bf j}= {\bf v}_{gf},\\
  &\quad {\bf j} = {\bf G}_{gffg}^{-1} {\bf Z}_{gffg} {\bf G}_{gffg}^{-1} \tilde{\bf j},
\end{align}
where
\begin{align*}
  &{\bf Z}_{gffg} = \sum_{i=1}^2
  \begin{bmatrix}
    {\bf K}_{i}^{gf}   & -\eta_{i} {\bf T}_{i}^{gg} \\
    \frac{1}{\eta_{i}}{\bf T}_{i}^{ff} & {\bf K}_{i}^{fg}
  \end{bmatrix},\\
  &{\bf G}_{gffg} = \\
  &\begin{bmatrix}
    \int_\Gamma\left(\hat{\Field{n}}\times \Field{g}_i\left(\Field{r}\right)\right)\cdot \Field{f}_j\left(\Field{r}\right) d\Field{r} & 0 \\
    0 &  \int_\Gamma\left(\hat{\Field{n}}\times \Field{f}_i\left(\Field{r}\right)\right)\cdot \Field{g}_j\left(\Field{r}\right) d\Field{r}
  \end{bmatrix},
\end{align*}
One of advantages of this alternative in the formulation of \eqref{eq_precond_pmchwt_mixed_matrix} is that the matrix ${\bf Z}_{gffg} {\bf G}_{gffg}^{-1}$ is also expected to be well-conditioned since $({\bf Z}_{gffg} {\bf G}_{gffg}^{-1})^2$ is well-conditioned as one sees from the Calder\'{o}n's identity and  equation \eqref{eq_precond_pmchwt_mixed_matrix}\cite{niino-jcp}.
This observation leads to the following well-conditioned discretized integral equation
\begin{align}\label{eq_precond_pmchwt_gram}
  &{\bf Z}_{gffg} {\bf G}_{gffg}^{-1}\tilde{\bf j} = {\bf v},\\
  \nonumber &\quad {\bf j} = {\bf G}_{gffg}^{-1}\tilde{\bf j}.
\end{align}
In section \ref{sec_orthogonalization}, we propose a CBFM combined with the CMP, which is based on the preconditioner in \eqref{eq_precond_pmchwt_gram}.

\section{CBFM with Calder\'{o}n Preconditioning}\label{sec_cbfm_with_cmp}
In this section we formulate the CBFM and its acceleration with the Calder\'{o}n preconditioning.
\subsection{CBF Generation} \label{sec_CBF_Generation}
\begin{figure}[!t]
  \centering
      \includegraphics[keepaspectratio, scale=0.4]{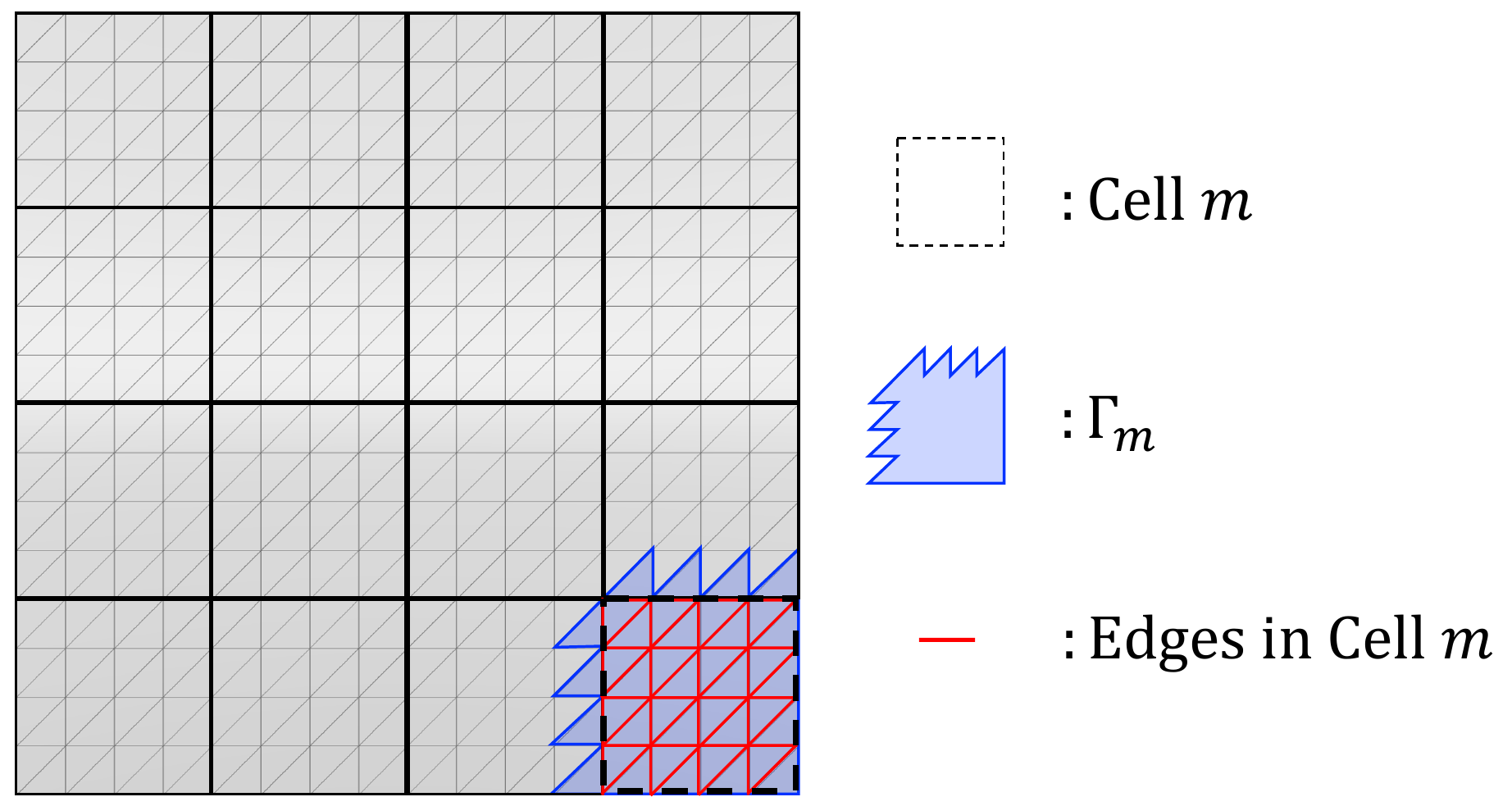}
  \caption{Division of scatterer into cells. }
  \label{fig_cell_division}
\end{figure}

In the CBFM \cite{Prakash2003}--\cite{Maaskant2008}, the scatterer is divided into $N^{\rm Cell}$ cells and CBFs are generated for each cell. 
\figuref\ref{fig_cell_division} shows an example of cell division of a scatterer. 
Thick black lines in the left part of the figure represent the cell boundaries for $N^{\rm Cell}=16$. 
Red edges in the right part, together with the associated triangles, indicate those contained in the cell enclosed by dashed lines in the left part. 
In addition, the region containing the support for all RWGs in cell $m$, indicated by the blue region, is defined as $\Gamma_{m}$. 
Note that an edge always belongs to one and only one cell.

The CBFs defined in cell $m$ are represented by linear combinations of the RWGs contained in that cell. 
The $n$th CBF $\Field{c}_{mn}\left(\Field{r}\right)$ in cell $m$ is expressed by 
\begin{align}
  \label{eq_cbf_JM}
  \Field{c}_{mn}(\Field{r})=
    \sum_{i=1}^{N_{m}}c_{mni}\Field{f}_{\Lambda_{mi}}\left(\Field{r}\right) \quad (n=1,\cdots L_m), 
\end{align}
where $N_m$  and $\Lambda_{mi}$ are the toal number and index of the RWG functions in cell $m$,
$L_m$ is the total number of CBFs, 
and $c_{mni}$ is the coefficient of the RWG for the CBF $\Field{c}_{mn}(\Field{r})$, respectively.
The matrix of the coefficients $c_{mni}$ in cell $m$ is denoted by ${\bf C}_{m}\CSize{N_{m}}{L_{m}}$;
\begin{align*}
  \begin{bmatrix}
    {\bf C}_{m}
  \end{bmatrix}_{in} = c_{mni}.
\end{align*}
In this paper, the coefficient ${\bf C}_{m}$ is also referred to as ``CBF'' unless it leads to confusion.

One can obtain CBFs by solving scattering problems with multiple incident fields.
In this paper, plane waves propagating in $s$ directions sampled at an appropriate interval on the unit sphere are used as incident fields. 
For instance, if $N_{\theta}$ ($N_{\phi}$) samples are used in the $\theta$ ($\phi$) direction and $N_{\rm p} (=1, 2)$ orthogonal polarization samples are considered, we have $s=N_{\theta}N_{\phi}N_{\rm p}$. 
From the $s$ plane waves, we compute the sets of coefficients ${\bf J}_{m}^{\rm J}$ and ${\bf J}_{m}^{\rm M}\CSize{N_m}{s}$, which correspond to solutions of scattering problems with the electric and magnetic incident fields illuminating the cell $m$.
We then obtain the CBFs ${\bf C}_{m}^{\rm J, M}$ by orthogonalizing ${\bf J}_{m}^{\rm J, M}$.
In the rest of this subsection we discuss how to determine ${\bf J}_{m}^{\rm J, M}$. The orthogonalization process to obtain ${\bf C}_{m}^{\rm J, M}$ from ${\bf J}_{m}^{\rm J, M}$ is described in section \ref{sec_orthogonalization}.

There are several options for generating ${\bf J}_{m}^{\rm J, M}$, depending on whether or not mutual couplings between cells is considered. 
The primary CBFs are the macro basis functions generated without taking into account the mutual couplings between cells.
The expansion coefficients ${\bf J}_{m}^{\rm M}$ and ${\bf J}_{m}^{\rm J}$ for the primary CBF can be calculated from the following equation \cite{Lucente2008}:
\begin{align}
  \label{eq_pcbf}
  {\bf Z}_{mm}{\bf J}_{m}&={\bf V}_{m} \\
  {\bf Z}_{mm} &= \sum_{i=1}^{2}
  \begin{bmatrix}
    \begin{pmatrix}{\bf K}_{i}^{ff}\end{pmatrix}_{mm} & -\eta_{i}\begin{pmatrix} {\bf T}_{i}^{ff}\end{pmatrix}_{mm} \\
    \frac{1}{\eta_{i}}\begin{pmatrix}{\bf T}_{i}^{ff}\end{pmatrix}_{mm} & \begin{pmatrix}{\bf K}_{i}^{ff}\end{pmatrix}_{mm} \\
  \end{bmatrix} \nonumber \\
  {\bf J}_{m} &=
  \begin{bmatrix}
    {\bf J}_{m}^{\rm M} \\
    {\bf J}_{m}^{\rm J} \\
  \end{bmatrix}\nonumber \\
  {\bf V}_{m} &=
  \begin{bmatrix}
    {\bf V}_{m}^{\rm E} \\
    {\bf V}_{m}^{\rm H} \\
  \end{bmatrix}\nonumber 
\end{align}
where $\begin{pmatrix} {\bf K}_{i}^{ff} \end{pmatrix}_{mm} \CSize{N_{m}}{N_{m}}$ denotes the submatrix of ${\bf K}_{i}^{ff}$ formed by $N_{m}$ RWGs in cell $m$, etc., and ${\bf V}_{m}^{\rm E, H}\CSize{N_{m}}{s}$ is the matrix of incident plane waves ${\bf v}$ for cell $m$. 
The primary CBFs are particularly useful for analyzing scatterers with periodic or quasi-periodic structures, such as array antennas or frequency selective surfaces whose components are separated from each other.
In such cases, it is convenient to configure each cell to contain one or several scatterers.

When one has to divide a single scatterer into several cells, however, equation \eqref{eq_pcbf} for each cell corresponds to integral equations defined on a divided open surface, which is not a valid PMCHWT formulation.
For such cases, one can use the IPCBF, which one may interpret as solutions of the whole system obtained by a block iterative solver.

The IPCBF \cite{Tanaka2016}, \cite{Tanaka2017} is one of CBFs which take into account the coupling between cells.
It is known that a reasonable RCS can be calculated from a small number of IPCBFs if $s$ incident fields are properly selected from the coordinate plane for the RCS pattern calculation \cite{Tanaka2019}.  
This number $s$ can be taken much smaller than the number $S$ of the directions for calculating the RCS as will be shown in section \ref{sec_num_ex}. 
The coefficient vector set for the IPCBF generation with $p$ iterations consisting of submatrices ${\bf J}_{m}^{{\rm M}, (p)} \CSize{N_{m}}{s}$ and ${\bf J}_{m}^{{\rm J}, (p)}{s} \CSize{N_{m}}{s}$ corresponding to the magnetic and electric currents, defined in \cite{Tanaka2016,Tanaka2017} and denoted by ${\bf J}_{m}^{(p)} \CSize{2N_{m}}{s}$, can be obtained from the following equation;
\begin{align}
  \label{eq_ipcbf}
  {\bf J}_{m}^{(p)}&=
  \begin{bmatrix}
    {\bf J}_{m}^{{\rm M}, (p)} \\
    {\bf J}_{m}^{{\rm J}, (p)} \\
  \end{bmatrix} \nonumber \\
  &={\bf Z}_{mm}^{-1}\left({\bf V}_{m}-\sum_{\substack{n=1\\n\neq m}}^{M}{\bf Z}_{mn}{\bf J}_{n}^{(p-1)} \right), \\
  {\bf Z}_{mn} &= \sum_{i=1}^{2}
  \begin{bmatrix}
    \begin{pmatrix}{\bf K}_{i}^{ff}\end{pmatrix}_{mn} & -\eta_{i}\begin{pmatrix} {\bf T}_{i}^{ff}\end{pmatrix}_{mn} \\
    \frac{1}{\eta_{i}}\begin{pmatrix}{\bf T}_{i}^{ff}\end{pmatrix}_{mn} & \begin{pmatrix}{\bf K}_{i}^{ff}\end{pmatrix}_{mn} \\
  \end{bmatrix} \nonumber 
\end{align}
where $p>0$ is a small integer, ${\bf Z}_{mn} \CSize{2N_{m}}{2N_{n}}$ is the submatrix of the impedance matrix $\bf Z$ representing the interaction between the elements in cells $m$ and $n$, and $N_{n}$ is the number of the RWGs in cell $n$. 
We note that \eqref{eq_ipcbf} gives the algorithm of the block Jacobi method. This shows that the IPCBF utilizes approximate iterative solutions of the MoM as basis functions. 
For the convergence of stationary methods such as the block Jacobi method, however, the impedance matrix must be diagonally dominant,
which may not always be the case in EM applications.
It is known, however, that some non-stationary methods such as the generalized minimal residual (GMRES) algorithm \cite{Saad2003} do not require such dominance for convergence. 
In this paper, we therefore propose to construct the coefficient sets ${\bf J}_m^{\rm J, M}$ for IPCBFs by solving each of plane wave problems for the $s$ incident directions using GMRES with a relatively large allowance $\delta_{\rm r}$ for the residual norm.

Tanaka et al. \cite{Tanaka2019} have shown that the choices of the directions and the number of incident plane waves have a significant impact on the accuracy in the CBFM analysis using \eqref{eq_ipcbf}.
They have also shown that the accuracy of the final solution can be controlled by the magnitude of the residual norm $\delta_{\rm r}$ in the iterative calculation of coefficient vector in \eqref{eq_coef_vector}.
In section~\ref{sec_num_ex}, we discuss the impact of $\delta_{\rm r}$ on the accuracy of the solution and the computational time in the GMRES version of IPCBF.

\subsection{CBF orthogonalization considering the property of the electromagnetic currents} \label{sec_orthogonalization}
The coefficient sets ${\bf J}_m^{\rm J, M}$ obtained in the last subsection cannot be directly used as basis functions since they can be linearly-dependent.
Hence we orthogonalize CBFs in cell $m$ using SVD \cite{Lucente2008}--\cite{Maaskant2008}. 
In the conventional CBFMs, the SVD is usually applied to the set ${\bf J}_{m}^{\rm J}$ and ${\bf J}_{m}^{\rm M}$ separately. 
The orthogonalized set for the electric current in cell $m$, for example, is obtained with the SVD as follows:
\begin{align}
  \label{eq_gram_CBF_SVD}
  {\bf J}_{m}^{\rm J} 
  =
  \begin{bmatrix}
    {\bf U}_{1} & {\bf U}_{2}
  \end{bmatrix}
  \begin{bmatrix}
    {\bf \Sigma}_{1} &  \\
    & {\bf \Sigma}_{2}
  \end{bmatrix}
  \begin{bmatrix}
    {\bf V}_{1}^{H} \\
    {\bf V}_{2}^{H}
  \end{bmatrix},
\end{align}
where ${\bf U}_i$ and ${\bf V}_i$ $(i=1,2)$ are unitary matrices having the singular vectors as their columns
and $(\cdot)^{H}$ stands for the adjoint.
Also, $\Sigma_i$ is a diagonal matrix with singular values as its diagonal components in descending order.
$\Sigma_1$ contains all the singular values exceeding the pre-defined threshold while singular values in $\Sigma_2$ are less than the threshold.
We now take columns of ${\bf U}_{1}$ as linearly independent CBFs for representing the electric currents in the cell $m$.
${\bf J}_{m}^{\rm M} $ is calculated similarly.

We note that CBFs for electric and magnetic currents thus obtained are not related.
However, it is natural to ask 
if these basis functions can have a duality structure as do the RWG and BC basis functions, for example.
If this is possible, one may establish a method of analysis that combines the advantages of CBFM and CMP, with fewer unknowns and faster convergence than the conventional method.
In what follows we propose a method for generating CBFs in which a duality relationship is satisfied. 

We now start with the definition of the duality. 
The CBFs $\Field{c}_{mi}^{\rm M}$ ($i=1,\cdots, L_{m}^{\rm M}$) and $\Field{c}_{mj}^{\rm J}$  ($j=1,\cdots, L_{m}^{\rm J}$) for magnetic and electric currents within a cell $m$ are said to be dual if  $L_{m}^{\rm M}= L_{m}^{\rm J}$, and the matrix
\begin{align}
  \label{eq_duality}
  \int_{\Gamma_{m}} \left(\hat{\Field{n}}\times\Field{c}^{\rm M}_{mi}(\Field{r})\right)\cdot\Field{c}^{\rm J}_{mj}(\Field{r})d\Field{r} 
\end{align}
is diagonal with positive diagonals. 
To determine CBFs satisfying this duality, we compute the Gram matrix ${\bf G}_{m}^{\prime}$ defined by
\begin{align}
  \label{eq_gram_CBF}
  {\bf G}_{m}^{\prime} &= \left({\bf J}_{m}^{\rm M}\right)^{H} {\bf G}_{m} {\bf J}_{m}^{\rm J} \CSize{s}{s},
\end{align}
where ${\bf G}_{m} \CSize{ N_{m} }{ N_{m} }$ is the Gram matrix of the RWG functions $\Field{f}_{i}$ contained in cell $m$:
\begin{align}
  \label{eq_gram_RWG}
  \begin{bmatrix}
    {\bf G}_{m}
  \end{bmatrix}_{ij}
  &= \int_{\Gamma_{m}} \left(\hat{\Field{n}}\times\Field{f}_{\Lambda_{mi}}\left(\Field{r}\right)\right)\cdot\Field{f}_{\Lambda_{mj}}\left(\Field{r}\right)d\Field{r}.
\end{align}
We orthogonalize CBFs by applying the SVD to the Gram matrix ${\bf G}_{m}^{\prime}$ and cutting off singular vectors associated with small singular values as follows:
\begin{align}
  \label{eq_gram_CBF_SVD}
  {\bf G}_{m}^{\prime} &= \left({\bf J}_{m}^{\rm M}\right)^{H} {\bf G}_{m} {\bf J}_{m}^{\rm J}  \nonumber \\
  &=
  \begin{bmatrix}
    {\bf U}_{\rm L} & {\bf U}_{\rm S}
  \end{bmatrix}
  \begin{bmatrix}
    {\bf \Sigma}_{\rm L} &  \\
    & {\bf \Sigma}_{\rm S}
  \end{bmatrix}
  \begin{bmatrix}
    {\bf V}_{\rm L}^{H} \\
    {\bf V}_{\rm S}^{H}
  \end{bmatrix} \nonumber \\
  &= {\bf U}_{\rm L} {\bf \Sigma}_{\rm L} {\bf V}_{\rm L}^{H} + {\bf U}_{\rm S} {\bf \Sigma}_{\rm S} {\bf V}_{\rm S}^{H} \nonumber \\ 
  &\simeq {\bf U}_{\rm L} {\bf \Sigma}_{\rm L} {\bf V}_{\rm L}^{H} \CSize{s}{s} \left( {\bf \Sigma}_{\rm L} \CSize{L_{m}}{L_{m}}  \right),
\end{align}
where $L_{m}$ is the number of singular values whose ratio to the largest singular value is greater than or equal to the threshold value $\delta_{\rm SVD}$. 
The size of both of the matrices ${\bf U}_{\rm L}$ and ${\bf V}_{\rm L}$ is thus $s \times L_{m}$. 
The CBFs ${\bf C}_{m}^{\rm M}$ and ${\bf C}_{m}^{\rm J}$ corresponding to the electric and magnetic currents are now obtained by
\begin{align}
  \label{eq_CBF}
  {\bf C}_{m}^{\rm M} &= {\bf J}_{m}^{\rm M} {\bf U}_{\rm L} \CSize{N_{m}^{\rm M}}{L_{m}}\\ 
  {\bf C}_{m}^{\rm J} &= {\bf J}_{m}^{\rm J} {\bf V}_{\rm L} \CSize{N_{m}^{\rm J}}{L_{m}}.
\end{align}
It is easy to see that the CBFs thus obtained satisfy the duality in \eqref{eq_duality}. Indeed, the matrix in \eqref{eq_duality} is computed as
\begin{align}
\left({\bf C}_{m}^{\rm M}\right)^{H} {\bf G}_{m}^{\prime} {\bf C}_{m}^{\rm J} &= \left({\bf J}_{m}^{\rm M} {\bf U}_{\rm L}\right)^{H}{\bf G}_{m}{\bf J}_{m}^{\rm J} {\bf V}_{\rm L} \nonumber \\
&= \left({\bf U}_{\rm L}\right)^{H} \left({\bf J}_{m}^{\rm M}\right)^{H} {\bf G}_{m} {\bf J}_{m}^{\rm J} {\bf V}_{\rm L} \nonumber \\
\label{eq_CBF_sv} &= {\bf \Sigma}_{\rm L}
\end{align}
which is diagonal and with positive diagonals.

The number of CBFs in the cell $m$, denoted by $N_{m}^{\rm CBF}$,
is $L_{m}^{\rm M}+L_{m}^{\rm J}=2L_{m}$; therefore, the total number of CBFs, denoted by $N^{\rm CBF}$, is $2\sum_{m=1}^{N^{\rm Cell}}L_{m}$.
In this way we obtain CBFs with the use of coefficient set ${\bf C}_{m}^{\rm M}$ and ${\bf C}_{m}^{\rm J}$, which are mutually dual. 
\subsection{Matrix Equation for the CBFM and Calder\'{o}n Multiplicative Preconditioner}
With the dual CBFs obtained in the last section, we apply the CMP introduced in \eqref{eq_precond_pmchwt_gram}.
For simplicity of explanation, we assume that elements associated with the RWG functions of smaller indices are in cells with smaller indices in order.
This condition does not have to be satisfied in the actual implementation, though.
We consider two sparse matrices ${\bf C}^{\rm JM} \CSize{N}{N^{\rm CBF}}$ and ${\bf C}^{\rm MJ} \CSize{N}{N^{\rm CBF}}$ defined as follows:  
\begin{align}
  \label{eq_CBF_matrix_JM}
  {\bf C}^{\rm JM}
  &\equiv
  \begin{bmatrix}
  {\bf C}_{1}^{\rm J} & & & & & \\
   & \ddots & & & & \\
   & & {\bf C}_{N^{\rm Cell}}^{\rm J} & & & \\
   & & & & {\bf C}_{1}^{\rm M} & & \\
   & & & & & \ddots & \\
   & & & & & & {\bf C}_{N^{\rm Cell}}^{\rm M}
  \end{bmatrix}, \\
  \label{eq_CBF_matrix_MJ}
  {\bf C}^{\rm MJ} 
  &\equiv
  \begin{bmatrix}
  {\bf C}_{1}^{\rm M} & & & & & \\
   & \ddots & & & & \\
   & & {\bf C}_{N^{\rm Cell}}^{\rm M} & & & \\
   & & & & {\bf C}_{1}^{\rm J} & & \\
   & & & & & \ddots & \\
   & & & & & & {\bf C}_{N^{\rm Cell}}^{\rm J}
  \end{bmatrix}.
\end{align}
By using \eqref{eq_CBF_matrix_JM} and \eqref{eq_CBF_matrix_MJ} we obtain the right preconditioner ${\bf G}^{\rm CBF}$,
which corresponds to the preconditioner ${\bf G}_{gffg}$ in \eqref{eq_precond_pmchwt_gram} with the RWG and BC functions replaced with the dual CBFs, given by
\begin{align}
  \label{eq_CBF_Gram_matrix}
  {\bf G}^{\rm CBF} &= \left( {\bf C}^{\rm JM} \right)^{H}{\bf G}_{ffff} {\bf C}^{\rm MJ}  \CSize{N^{\rm CBF}}{N^{\rm CBF}},
\end{align}
where ${\bf G}_{ffff}$ is defined by \eqref{eq_Gffff}. 
Note that the matrix ${\bf G}^{\rm CBF}$ is diagonal if no connected boundary is divided into multiple cells in the CBFM.
This is because the Gram matrix ${\bf G}^{\rm CBF}$ is block diagonal and each block corresponds to the Gram matrix in \eqref{eq_gram_RWG}.
However, ${\bf G}^{\rm CBF}$ is not diagonal when a connected boundary is divided into multiple cells. We will return to this issue at the end of this section.

The matrix equation in CBFM with the use of right preconditioner ${\bf G}^{\rm CBF}$ can be written as
\begin{align}
  \label{eq_CBFM_outer}
  \left[ \left( {\bf C}^{\rm JM} \right) ^{H} {\bf Z}_{ffff} {\bf C}^{\rm MJ} \right] \left({\bf G}^{\rm CBF}\right)^{-1}{\bf y}^{\rm CBF} &= \left( {\bf C}^{\rm JM} \right) ^{H} {\bf v}_{ff}, \\
  \label{eq_CBFM_inner}
  {\bf G}^{\rm CBF}{\bf j}^{\rm CBF} &= {\bf y}^{\rm CBF}. 
\end{align}
As mentioned in section \ref{sec_pre-cald}, an appropriate Gram matrix as a preconditioner improves the convergence of the matrix equation with {\bf K} operator as the diagonal components.
The proposed CBFM expressed in \eqref{eq_CBFM_outer} and \eqref{eq_CBFM_inner} is expected to inherit this property.
Finally, the expansion coefficient vector ${\bf j}$ in \eqref{eq_coef_vector} is obtained from the solution ${\bf j}^{\rm CBF}$ as follows:
\begin{align}
  \label{eq_trans_current}
  {\bf j} = {\bf C}^{\rm MJ}{\bf j}^{\rm CBF}, 
\end{align}
and the electric and magnetic currents can be calculated with \eqref{eq_J_expand} and \eqref{eq_M_expand}. 
Note that the matrix products in \eqref{eq_CBFM_outer} and \eqref{eq_CBFM_inner} can be computed efficiently
since ${\bf C}^{\rm MJ}$, ${\bf C}^{\rm JM}$ are sparse, ${\bf G}^{\rm CBF}$ is diagonal and the matrix-vector product with ${\bf Z}_{ffff}$ can be accelerated by fast methods such as the FMM \cite{Chew2001}-\cite{Song1997}.
Consequently the proposed method is expected to have good properties in terms of both convergence and computational cost. 
From another perspective, the proposed method has an advantage that
one can easily implement it by adding a few components to an existing in-house program of the MoM. Indeed, the proposed method only needs to calculate ${\bf C}^{\rm MJ}$, ${\bf C}^{\rm JM}$ and ${\bf G}^{\rm CBF}$ anew and put them into the conventional iterative method as preconditioners, as is evident in \eqref{eq_CBFM_outer} and \eqref{eq_CBFM_inner}.

If a connected boundary is divided into multiple cells, the Gram matrix ${\bf G}^{\rm CBF}$ is no longer block diagonal since the Gram matrix has non-zero elements corresponding to interactions between two cells.
In this case the Gram matrix of the CBF, i.e., ${\bf G}^{\rm CBF}$ in \eqref{eq_CBF_Gram_matrix}, is not diagonal and could be ill-conditioned in the worst case.
However we consider that the preconditioning in \eqref{eq_CBFM_outer} and \eqref{eq_CBFM_inner} is effective in many problems
since the number of non-zero components outside the diagonal blocks is much smaller than the total DOF.
In addition, $\left({\bf G}^{\rm CBF}\right)^{-1}$ can be easily calculated with an iteration method using the diagonal preconditioner ${\bf D}^{\rm CBF}$ defined as
\begin{align}
  {\bf D}^{\rm CBF} &= 
  \begin{cases}
  {}
    1/\begin{bmatrix} {\bf G}^{\rm CBF} \end{bmatrix}_{ij}, & i=j \\
    0, & \mbox{otherwise}
  \end{cases}
  \CSize{N^{\rm CBF}}{N^{\rm CBF}}.
\end{align}
since ${\bf G}^{\rm CBF}$ is sparse and its diagonal elements, which are actually the singular values in \eqref{eq_CBF_sv}, are relatively large numbers.
We will verify numerically the above statement in section \ref{sec_Cylinder}.

\section{Numerical Results}\label{sec_num_ex}
In this section, we evaluate the validity of the proposed method by analyzing several scattering problems and comparing the results with those obtained with MoM. 
For a fair comparison with the proposed method, we use the PMCHWT formulation with the following commonly used operator arrangement in the MoM:
\begin{align}
  \label{eq_MoM_T}
  {\bf Z}_{\rm t}  {\bf D}^{-1} \tilde{\bf j}_{\rm t} &=  {\bf v}, \\
  {\bf j}_{\rm t} &= {\bf D}^{-1} \tilde{\bf j}_{\rm t}, \nonumber \\
  {\bf Z}_{\rm t}
  &=\sum_{i=1}^{2}
  \begin{bmatrix}
    -\eta_{i} {\bf T}_{i}^{gg} & {\bf K}_{i}^{gf}\\
    {\bf K}_{i}^{fg} & \frac{1}{\eta_{i}}{\bf T}_{i}^{ff} \\
  \end{bmatrix}, \nonumber \\
  {\bf j}_{\rm t} 
  &=
  \begin{bmatrix}
    \begin{bmatrix}
      \alpha^{\rm J}_{1} & \cdots & \alpha^{\rm J}_{N^{{\rm J}}}
    \end{bmatrix} & 
    \begin{bmatrix}
      \alpha^{\rm M}_{1} & \cdots & \alpha^{\rm M}_{N^{{\rm M}}} 
    \end{bmatrix}
  \end{bmatrix}^{T}, \nonumber 
\end{align}
where ${\bf D}$ is the diagonal preconditioner. 
This arrangement of the operators appears to be more natural than the proposed method in \eqref{eq_Z}
since the matrix ${\bf T}_i$ has large diagonal elements while ${\bf K}_i$ corresponds to a discretized compact operator which has smaller diagonal elements. 
In the following examples, \eqref{eq_pcbf} is used to obtain currents for generating primary CBFs except in \ref{sec_Cylinder}, where \eqref{eq_MoM_T} is used for generating IPCBFs. 
Throughout this section, we use the FMM for the multiplication of the impedance matrix and the GMRES for solving linear equations. 
In the example of IPCBFs in \ref{sec_Cylinder} the $({\bf G}^{\rm CBF})^{-1}$ operation in \eqref{eq_CBFM_outer} is carried out
with the GMRES using the empirically determined error tolerance of $1.0\times10^{-5}$. 
When the number of unknowns in a single cell is small, one can solve equation \eqref{eq_pcbf} with the LU decomposition. 
It can be used to obtain solutions faster for multiple incident waves. 
From the reason, we use the LU decomposition for solving the equation \eqref{eq_pcbf} to obtain the primary CBFs. 
\subsection{Sphere Array}\label{sec_sphere_array}
\subsubsection{Single Sphere}
\begin{figure}[!t]
  \centering
    \includegraphics[keepaspectratio, scale=0.4]{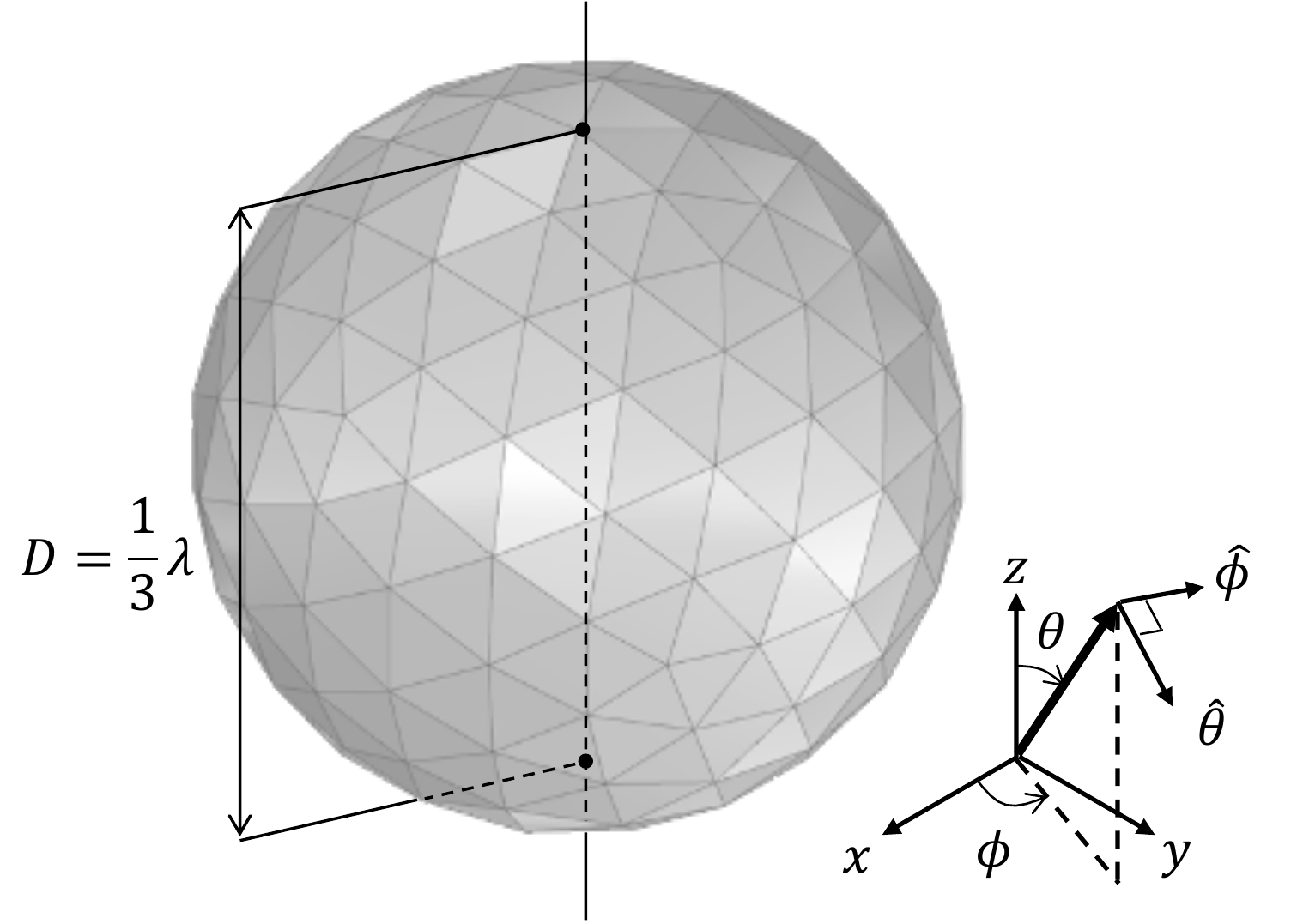}
  \caption{Sphere}
  \label{fig_Sphere_R5mm_Tess2}
\end{figure}

The first example is a sphere of diameter $\lambda/3$ as shown in \figuref\ref{fig_Sphere_R5mm_Tess2}, whose
relative permittivity $\epsilon_{\rm r}$ is varied from $1.0$ to $6.0$. 
The number of unknowns $N$ is $960$.  
The primary CBFs are used for the analysis. 

\tableref\ref{table_condition_sphere} shows parameters for the incident waves used for calculating CBFs.
We generate CBFs by solving \eqref{eq_pcbf} with the incident plane waves propagating in the directions given by the spherical coordinates
\begin{align}
 \label{eq_CBF_sph_coord} (\theta, \phi) &= (\theta_{\rm s}+n_\theta\Delta\theta, \phi_{\rm s}+n_\phi\Delta\phi),\\
 \nonumber  &\textrm{for}\quad n_\theta = 0,\ldots , N_\theta \quad  n_\phi=0, \ldots , N_\phi.
\end{align}
In this analysis, we consider one single cell which includes the entire sphere, thus not
dividing the scatterer into cells.
An analysis of this type may not be of much practical interest, but helps us to understand the characteristics of the CBFM, as we will see.
\setlength{\tabcolsep}{1.5mm}
\begin{table}[!t]
  \caption{Parameters for Calculating CBFs for Sphere Array}
  \centering
  \begin{tabular}{ccccccccc} \hline
    $\theta_{\rm s}$  & $\Delta\theta$ & $N_{\theta}$ & $\phi_{\rm s}$ & $\Delta\phi$ & $N_{\phi}$ & Pol.                                & Cell ($\lambda$) \\ \hline \hline
    $0^{\circ}$          & $30^{\circ}$     & $12$               & $0^{\circ}$      & $30^{\circ}$  & $6$            & $\hat{\theta},\hat{\phi}$  & $1/2$                \\ \hline
  \end{tabular}
  \label{table_condition_sphere}
\end{table}

\begin{figure}[!t]
  \centering
    \includegraphics[keepaspectratio, scale=0.55]{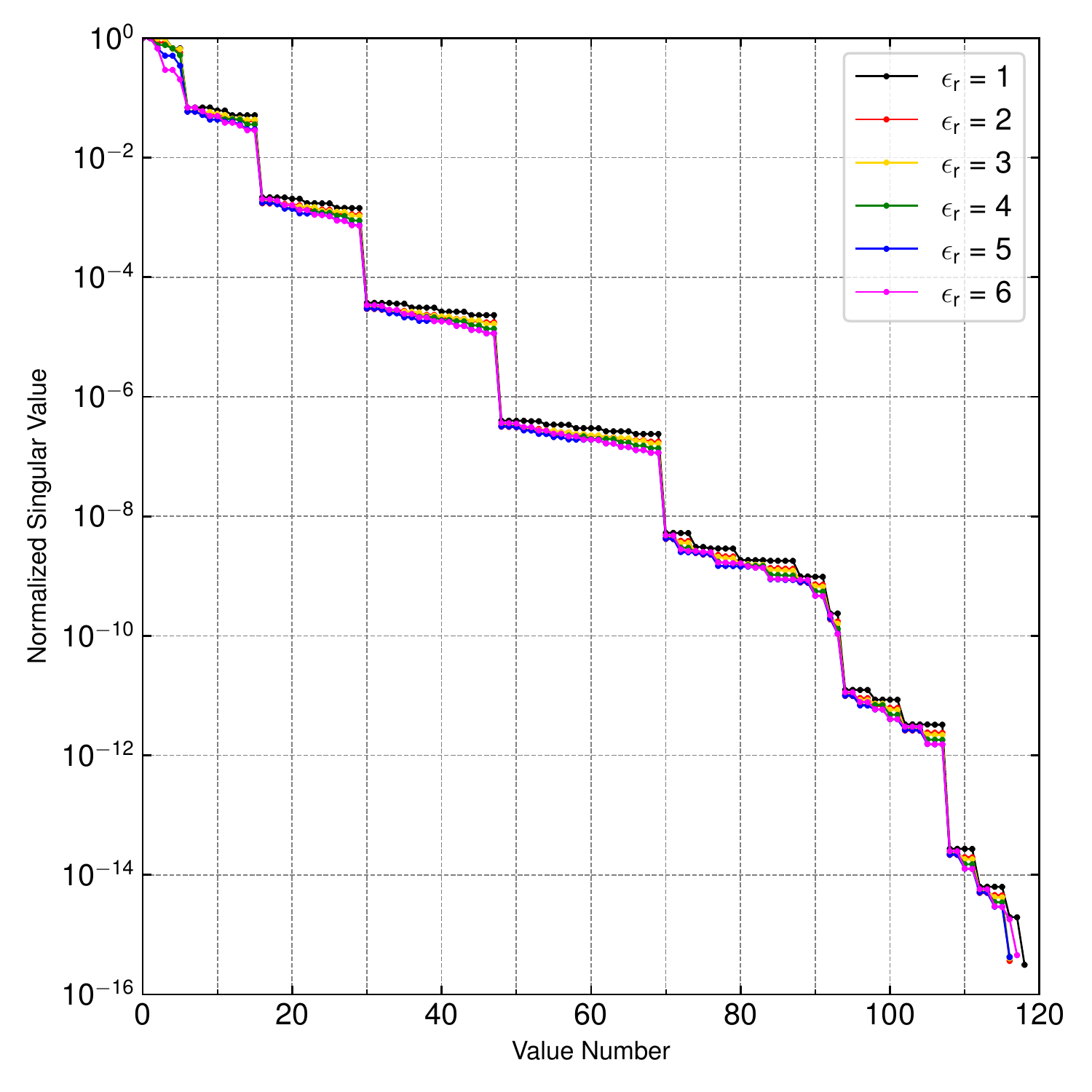}
  \caption{Normalized singular value of the sphere}
  \label{fig_svd_Sphere_R5mm_Tess2}
\end{figure}

\figuref\ref{fig_svd_Sphere_R5mm_Tess2} shows the distribution of singular values
in \eqref{eq_gram_CBF_SVD}.
In this figure, the horizontal axis shows the number of singular values assigned in descending order of the magnitudes, and the vertical axis shows the magnitude of each singular value normalized by their maximum value. 
We see that these singular values are divided into several groups regardless of the dielectric constant. 
Numbering these groups from left to right as $l (=1,2,\cdots)$, we see that the number of singular values belonging to the $l$th group is $2(2 l+1)$ which is the dimension of the vector spherical harmonics of the $l$th-order.
Motivated by this observation, we first consider up to $70$ singular values belonging to the group of $l = 1,\cdots,5$. 
Since $70$ CBFs are defined for both electric and magnetic currents, we have $N^{\rm CBF}=140$. 

\begin{figure}[!t]
  \begin{tabular}{c}
    \begin{minipage}[t]{0.45\hsize}
      \centering
      \includegraphics[keepaspectratio, scale=0.27]{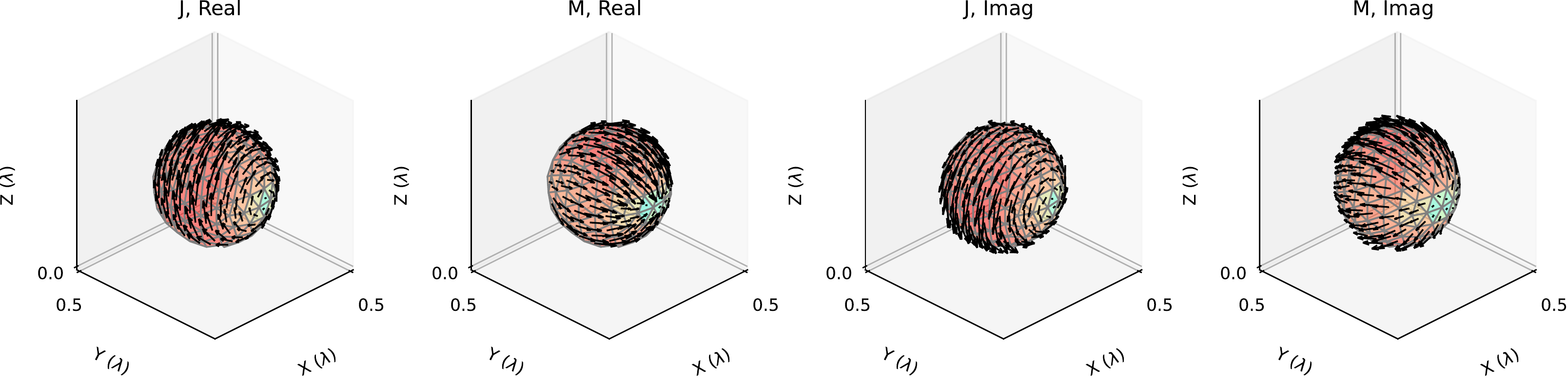}
      \subcaption{CBF for singular value \#$1$}
      \label{fig_Sphere_R5mm_Tess2_current0_CBF0_az45_el30}
    \end{minipage} \\
    \begin{minipage}[t]{0.45\hsize}
      \centering
      \includegraphics[keepaspectratio, scale=0.27]{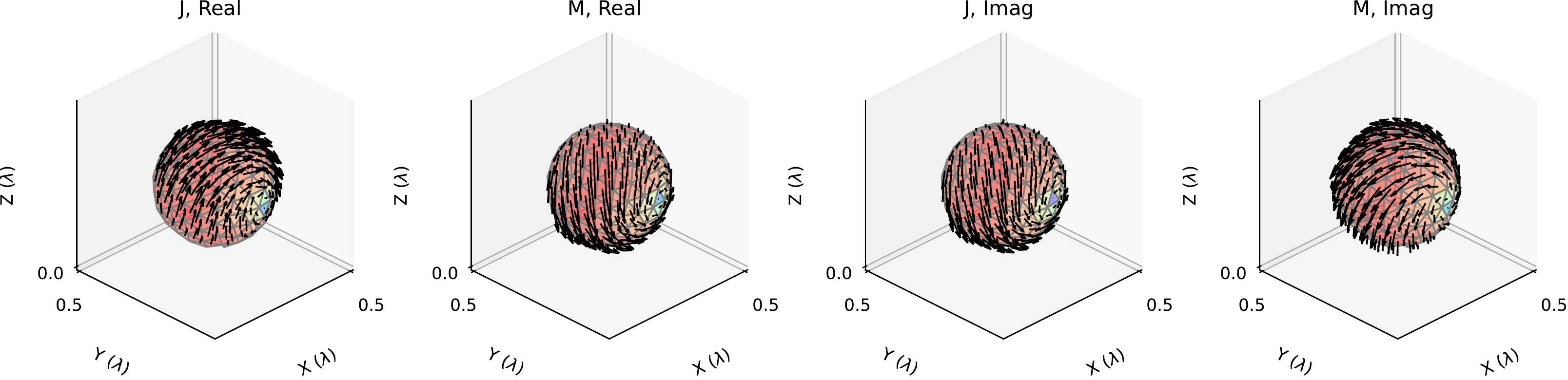}
      \subcaption{CBF for singular value \#$2$}
      \label{fig_Sphere_R5mm_Tess2_current0_CBF1_az45_el30}
    \end{minipage} \\
    \begin{minipage}[t]{0.45\hsize}
      \centering
      \includegraphics[keepaspectratio, scale=0.27]{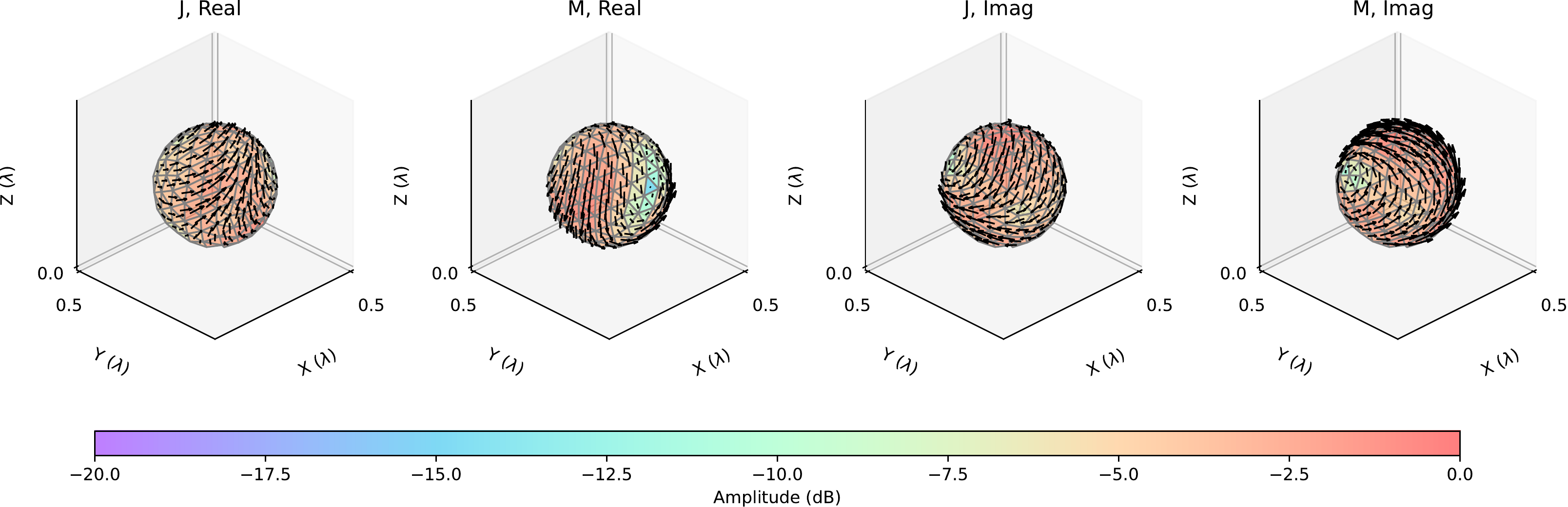}
      \subcaption{CBF for singular value \#$70$}
      \label{fig_Sphere_R5mm_Tess2_current0_CBF69_az45_el30}
    \end{minipage} 
  \end{tabular}
  \caption{CBF distributions on a sphere for $\epsilon_{\rm r}=1$. Black arrows in each figure represent current directions at the center of each mesh. \#$i$ corresponds to the singular value number. The two figures on the left show the real parts of ${\bf J}_{1}$ and ${\bf M}_{1}$, and the two right figures show their imaginary parts.}
  \label{fig_Sphere_R5mm_Tess2_current0}
\end{figure}

The proposed method is intended to orthogonalize the electric and magnetic CBFs as described in section \ref{sec_orthogonalization}. 
To check this, we compute CBFs at the barycentric coordinate of each triangle in a cell represented by \eqref{eq_cbf_JM}. 
\figuref\ref{fig_Sphere_R5mm_Tess2_current0} shows the real and imaginary parts of CBFs for electric and magnetic currents. 
\figuref\ref{fig_Sphere_R5mm_Tess2_current0_CBF0_az45_el30} (\figuref\ref{fig_Sphere_R5mm_Tess2_current0_CBF69_az45_el30}) gives the CBF distributions corresponding to the largest (smallest) singular value for $\epsilon_{\rm r}=1$. 
It is seen, indeed, that the electric and the corresponding magnetic CBFs are mutually orthogonal even up to the $70$th ones. 

\begin{figure}[!t]
  \centering
    \includegraphics[keepaspectratio, scale=0.55]{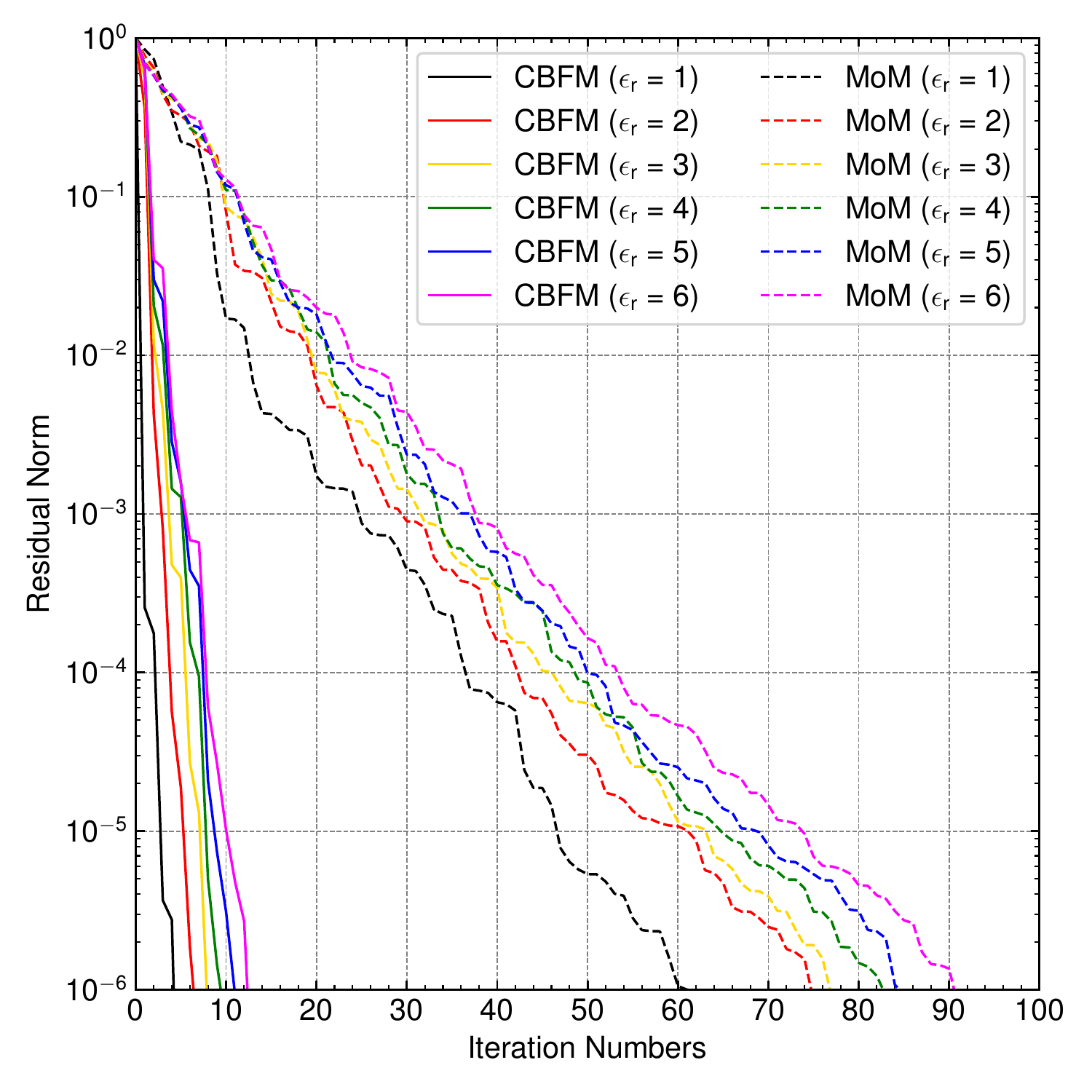}
  \caption{Convergence of the proposed CBFM with the CMP and the conventional MoM for a sphere. The incident field propagates in $-z$ direction}
  \label{fig_norm_Sphere_R5mm_Tess2}
\end{figure}

The convergence of the proposed CBFM with the CMP and the MoM using \eqref{eq_MoM_T} is shown in \figuref\ref{fig_norm_Sphere_R5mm_Tess2}, where the incident field propagates in $-z$ direction.
It can be clearly seen that the convergence of the proposed method is much improved compared to MoM. 
\subsubsection{Sphere Array}
\begin{figure}[!t]
  \centering
    \includegraphics[keepaspectratio, scale=0.4]{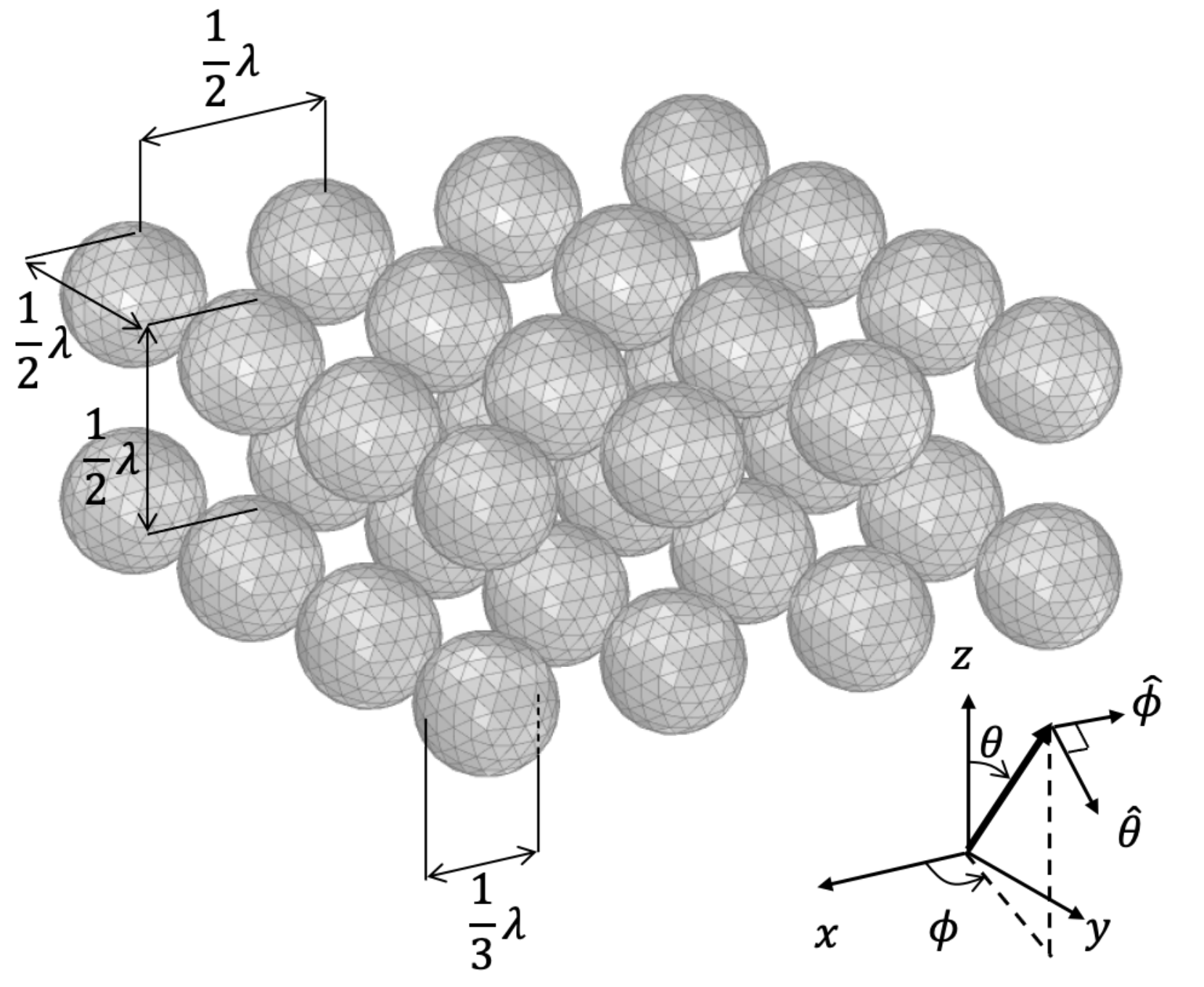}
  \caption{$4\times4\times2$ sphere array}
  \label{fig_Sphere_R5mm_Tess2_Array_4x4x2_15mm_space}
\end{figure}

Next, we analyze an array-shaped scatterer composed of spheres of the above mentioned shape aligned $4\times4\times2$ in the $x$, $y$, and $z$ directions as shown in \figuref\ref{fig_Sphere_R5mm_Tess2_Array_4x4x2_15mm_space}. 
The number of the unknowns $N$ is $30720$.  
The conditions of the incident field and the number of singular values considered for the generation of CBFs are the same as those in the case of a single sphere. 
In the previous analysis for a single sphere, the number of CBFs per cell was $140$; therefore $N^{\rm CBF}$ is $4480$ in this analysis. 
The number of unknowns will be reduced by a factor of about $6.9$ compared to the MoM. 

We first compare the convergence of the outer GMRES for three methods, i.e., the proposed method in \eqref{eq_CBFM_outer}, the proposed CBFM without using ${\bf G}^{\rm CBF}$, and the CBFM with the matrix arranged as in \eqref{eq_MoM_T}, i.e.
\begin{align}
  \label{eq_cbfm_Tdiag}
  \left[ \left( {\bf C}^{\rm JM} \right) ^{H} {\bf Z}_{\rm t} {\bf C}^{\rm JM} \right] {\bf j}_{\rm t}^{\rm CBF} &= \left( {\bf C}^{\rm JM} \right) ^{H} {\bf v} \\
  {\bf j}_{\rm t} &= {\bf C}^{\rm JM}{\bf j}_{\rm t}^{\rm CBF}. \nonumber
\end{align} 
For solving \eqref{eq_cbfm_Tdiag}, we do not use the right preconditioner ${\bf G}^{\rm CBF}$ in \eqref{eq_CBFM_outer} since this preconditioner is effective only when the ${\bf K}_i$-components are arranged diagonally as shown in \eqref{eq_precond_pmchwt_gram}. 

\begin{figure}[!t]
    \includegraphics[keepaspectratio, scale=0.55]{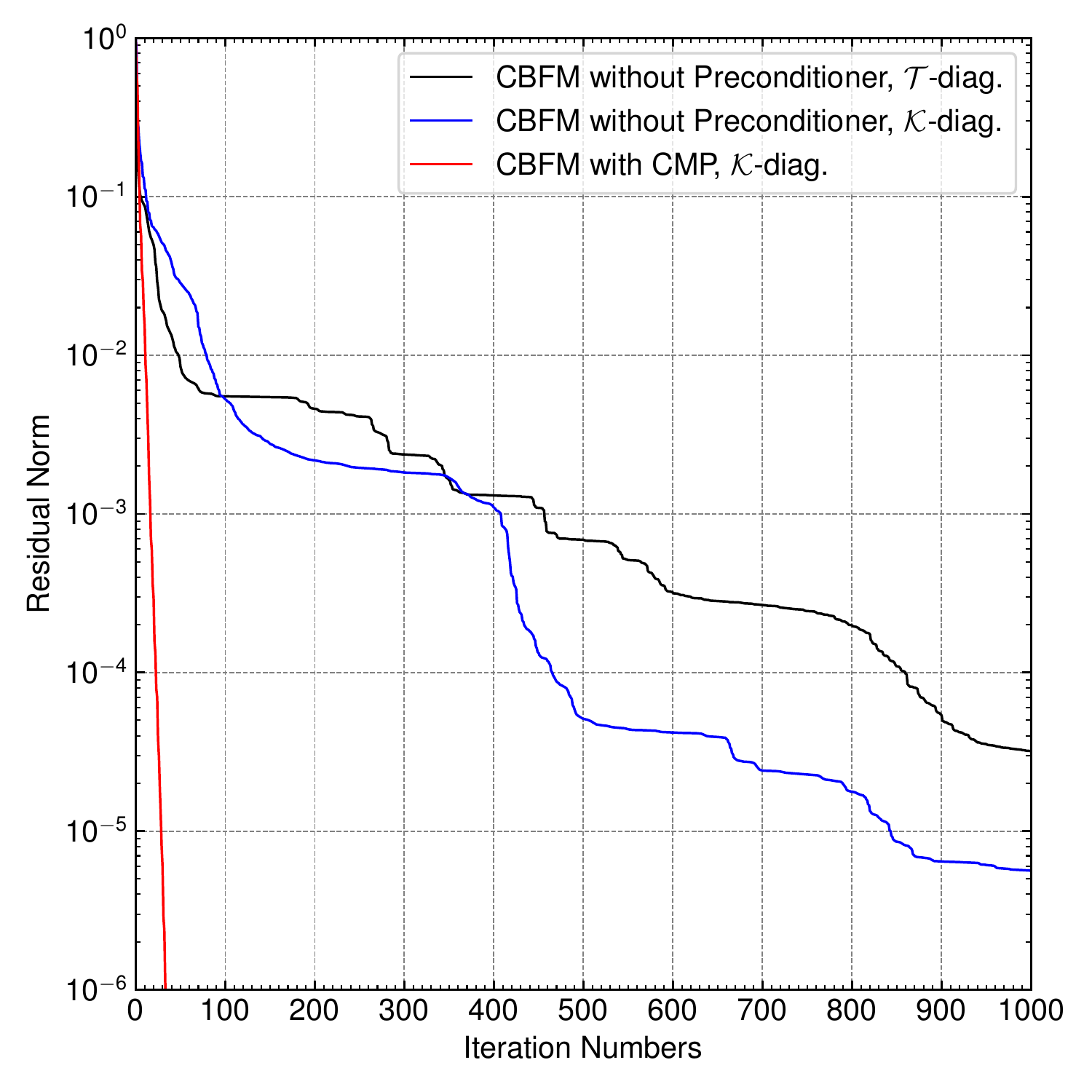}
  \caption{Convergence of the $4\times4\times2$ sphere array analysis for the incident wave with $\hat{\theta}$--polarization and $-z$ propagation direction. The relative permittivity of the scatterer is $\epsilon_{\rm r}=3$. }
  \label{fig_norm_Sphere_R5mm_Tess2_Array_4x4x2_15mm_space}
\end{figure}

\figuref\ref{fig_norm_Sphere_R5mm_Tess2_Array_4x4x2_15mm_space} shows the convergence of each method under the conditions that the incident plane wave propagates in $-z$ direction with $\hat{\theta}$-polarization and $\epsilon_{\rm r}=3$.
In the analyses without preconditioners, the residual norm $\delta_{\rm R}$ did not reach $10^{-6}$ even after $1000$ iterations. 
On the other hand, the proposed method with ${\bf G}^{\rm CBF}$  clearly improved the convergence since $\delta_{\rm R}$ reached $10^{-6}$ in less than $40$ iterations. 

\begin{figure}[!t]
    \includegraphics[keepaspectratio, scale=0.55]{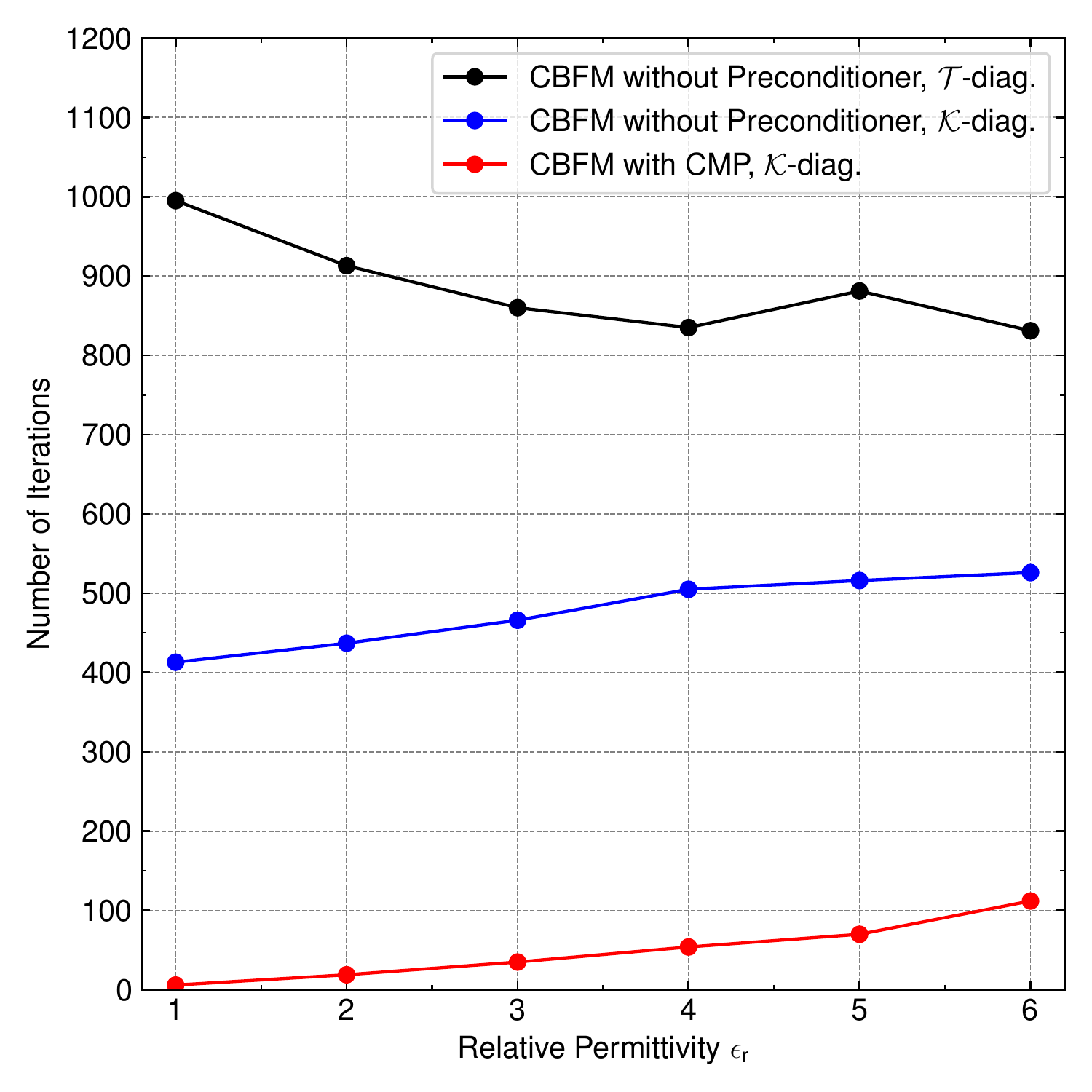}
  \caption{Relationship between the relative permittivity and the convergence of the scatterer}
  \label{fig_convergence_Sphere_R5mm_Tess2_Array_4x4x2_15mm_space_permittivity}
\end{figure}

Next, we check the convergence in the same scattering problem for variable relative permittivity $\epsilon_{\rm r}$. 
The threshold of the residual norm $\delta_{\rm R}$ is set to be $10^{-4}$, which the two methods without CMP could reach in \figuref\ref{fig_norm_Sphere_R5mm_Tess2_Array_4x4x2_15mm_space}.
\figuref\ref{fig_convergence_Sphere_R5mm_Tess2_Array_4x4x2_15mm_space_permittivity} shows the convergence properties of the three methods. 
The number of iterations for the proposed method using CMP is much smaller than those for other methods without CMP, although it gradually increases as the relative permittivity increases.
Indeed, the proposed approach could reduce the number  of iterations
to $1/5$ ($\epsilon_{\rm r}=6$) $\sim$ $1/70$ ($\epsilon_{\rm r}=1$) compared to the corresponding ${\bf K}$-diagonal CBFM without CMP.
Among methods without CMP, the ${\bf K}$-diagonal version seems to be  better than the ${\bf T}$-diagonal version as far as the convergence is concerned.

\begin{figure}[!t]
  \centering
    \includegraphics[keepaspectratio, scale=0.55]{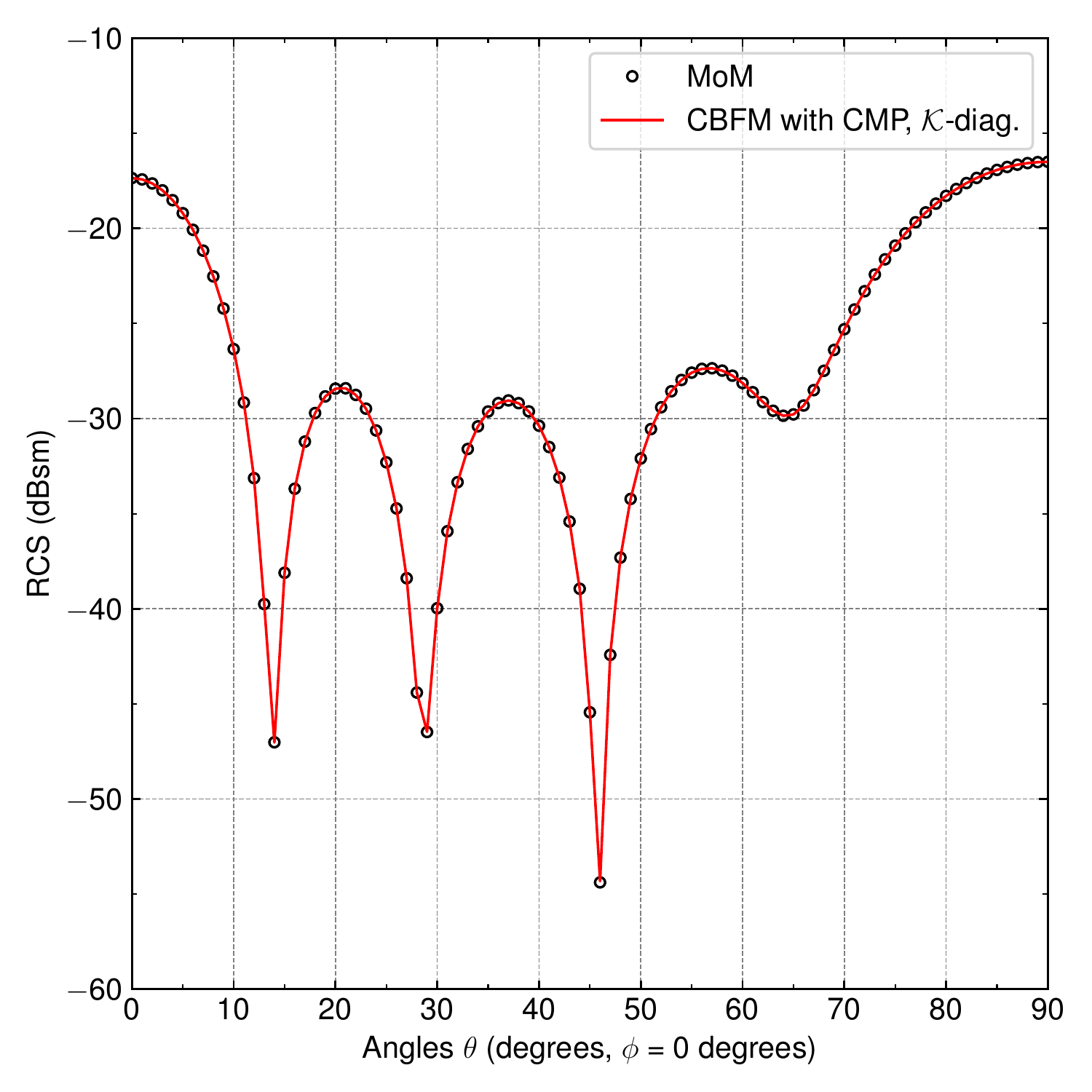}
  \caption{RCS of the sphere array ($4\times4\times2$) at $\epsilon_{\rm r}=3$}
  \label{fig_plot_rcs_Sphere_R5mm_Tess2_Array_4x4x2_15mm_space}
\end{figure}

\figuref\ref{fig_plot_rcs_Sphere_R5mm_Tess2_Array_4x4x2_15mm_space} compares the RCS patterns obtained with the proposed method and MoM, respectively, for $\epsilon_{\rm r}=3$. These results agree quite well. 
Indeed, the root mean square error~(RMSE) of the CBFM relative to the MoM defined by
\begin{align}
  \label{eq_RMSE}
  {\rm RMSE} = 10.*\log10\left(\frac{\sqrt{\frac{1}{N_{\theta}}\sum_{i=1}^{N_{\theta}}\left(\sigma_{i}^{\rm c}-\sigma_{i}^{\rm m}\right)^2}}{\underset{i}{\max}\{\sigma_{i}^{\rm m}\}-\underset{i}{\min}\{\sigma_{i}^{\rm m}\}}\right)\ \ {\rm dB} \nonumber
\end{align}
is $-39.98$ dB,
where $\sigma_i^{\rm c}$ and $\sigma_i^{\rm m}$ indicate the RCSs obtained with the CBFM and MoM in the $i$--th calculation, respectively. 
\subsubsection{Large Sphere Array}
As an example of a large-scale scattering analysis, we consider a sphere array consisting of $256$ spheres ($8 \times 16 \times 2$ in the $x$, $y$, and $z$ directions). 
The number of unknowns $N$ is $960 \times 256 = 245760$. 
In this analysis, we use the same condition and the relative permittivity ($\epsilon_{\rm r}=3$) as in the $4\times4\times2$ array case, hence we have $N^{\rm CBF}=35840$. 

\begin{figure}[!t]
    \includegraphics[keepaspectratio, scale=0.55]{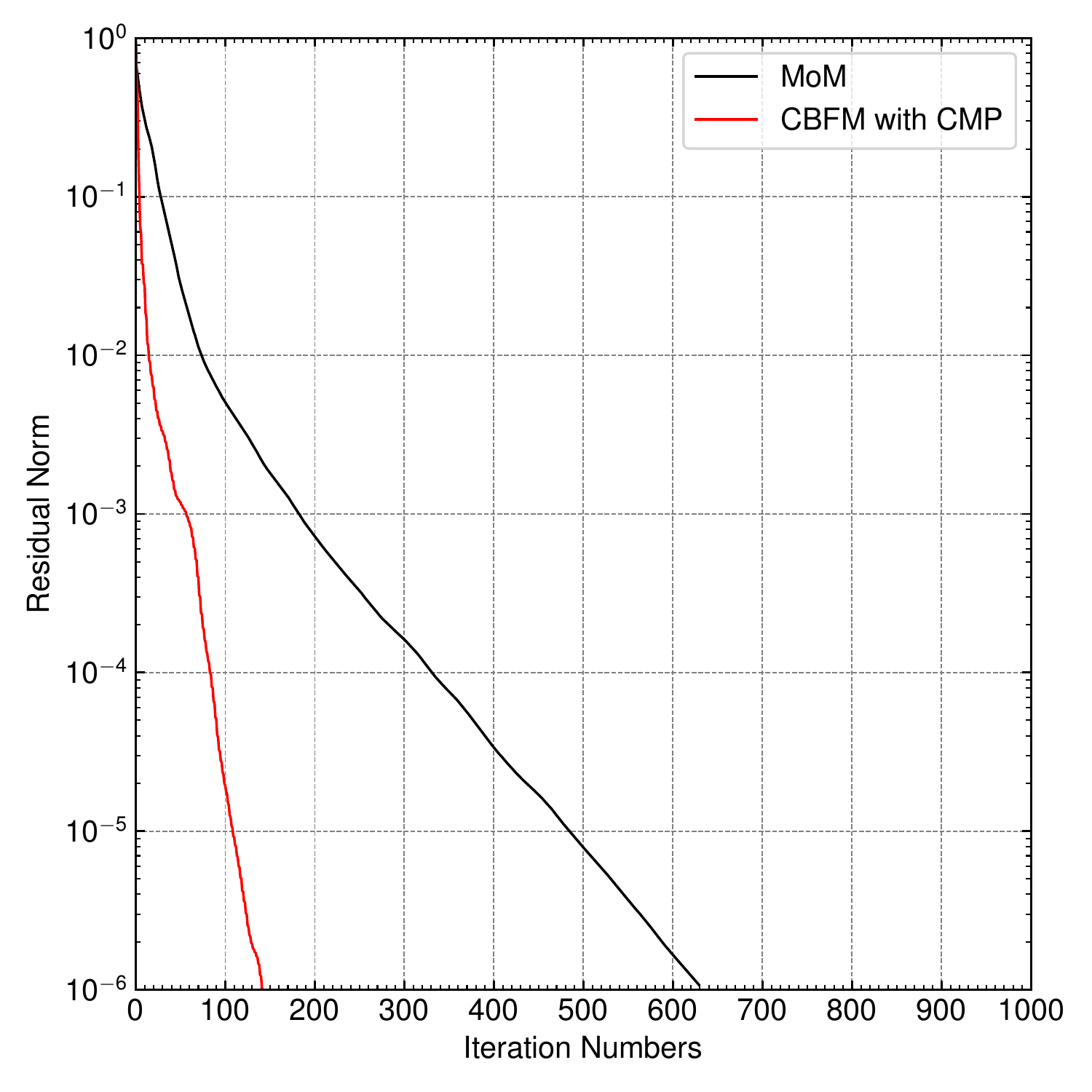}
  \caption{Convergence of the large sphere array analysis. }
  \label{fig_norm_Sphere_R5mm_Tess2_Array_8x16x2_15mm_space}
\end{figure}

\figuref\ref{fig_norm_Sphere_R5mm_Tess2_Array_8x16x2_15mm_space} shows the $\delta_{\rm R}$ (residual norm) vs the iteration number curves for the proposed method and MoM when the incident field is with the $\theta$-polarization and propagates in $-z$ direction.
In the proposed method, the number of iterations is $143$ when the residual norm $\delta_{\rm R}$ is $10^{-6}$, while it is about $634$ in the case of the MoM, which is about 4.4 times that of the proposed method.

\begin{figure}[!t]
  \centering
    \includegraphics[keepaspectratio, scale=0.55]{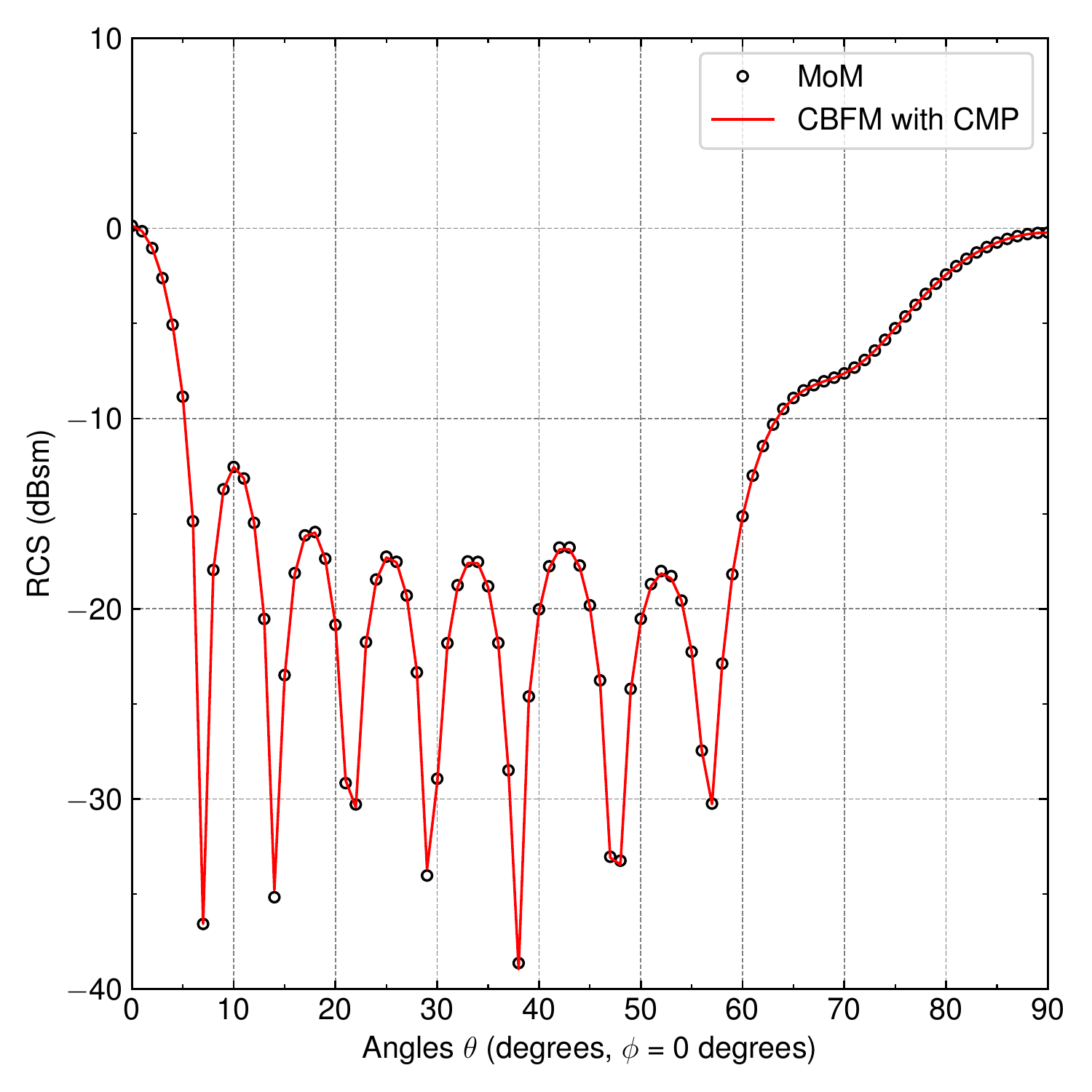}
  \caption{RCS of the large sphere array.}
  \label{fig_plot_rcs_Sphere_R5mm_Tess2_Array_8x16x2_15mm_space}
\end{figure}

\figuref\ref{fig_plot_rcs_Sphere_R5mm_Tess2_Array_8x16x2_15mm_space} shows the RCS pattern obtained using $\delta_{\rm R}<10^{-4}$ as the criterion of the convergence. 
The angle range and polarization to calculate the RCS pattern are $0^{\circ} \leq \theta \leq 30^{\circ}$, $\Delta\theta=1^{\circ}$, $\phi=0^{\circ}$ and $\hat{\theta}$. 
The RMSE is $-35.65$dB; hence the RCS obtained with the proposed CBFM is in good agreement with the MoM result. 
The computational times for generating the primary CBFs and final iteration, relative to the total computational time for the MoM for $91$ directions, are less than $0.02$ and approximately $0.21$, respectively. 
Hence the proposed method is more than four times faster than the conventional method in this analysis. 
These results show that the proposed method is effective particularly for large scale scatterers. 
\subsection{Cube Array} \label{sec_cube_array}
\begin{figure}[!t]
  \centering
    \includegraphics[keepaspectratio, scale=0.4]{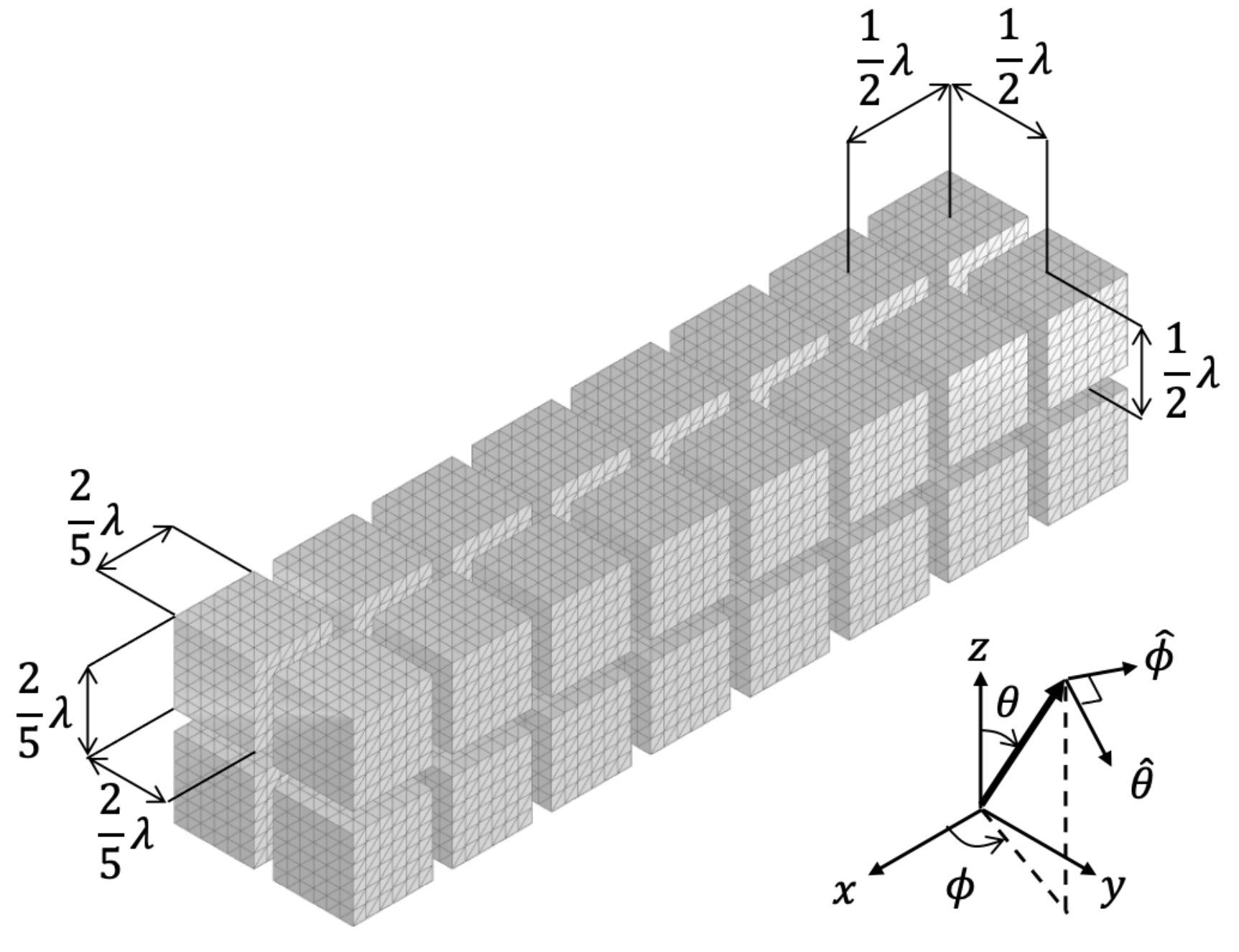}
  \caption{Cube array ($8\times2\times2$)}
  \label{fig_Cube_12_12_12_2mm_Array_8x2x2_15mm_space}
\end{figure}

The next analysis object is a cube array as shown in \figuref\ref{fig_Cube_12_12_12_2mm_Array_8x2x2_15mm_space}. 
The cubes are aligned $8\times2\times2$ in the $x$, $y$, and $z$ directions, respectively.
The number of the unknowns $N$ is $41472$ and 
the relative permittivity is assumed to be $3$. 

As in the previous section, we use primary CBFs. 
The cell is set up as in the previous analysis, i.e., a cubic cell with a side length  of $\lambda/2$ containing a cubic scatterer ($M=32$). 
The number of singular values to be considered is determined from the group number $l$ as mentioned in the previous section.
In this analysis, we also consider the impacts of the angular intervals of the generating incident fields and the threshold of the singular value on the accuracy and convergence of CBFM. The CMP in \eqref{eq_CBFM_outer} is used for all analyses. 
The primary CBFs are generated with plane waves, which propagate in the directions given by the spherical coordinates in \eqref{eq_CBF_sph_coord} with the parameters in \tableref\ref{table_condition_Cube_12_12_12_4mm_Array_8x2x2_15mm_space}.
\begin{table}[!t]
  \caption{Condition of CBFs for Cube Array}
  \centering
  \begin{tabular}{ccccccccc} \hline
    No.  & $\theta_{\rm s}$  & $\Delta\theta$ & $N_{\theta}$ & $\phi_{\rm s}$ & $\Delta\phi$ & $N_{\phi}$ & Pol.  & Group number $l$ \\ \hline \hline
    1 & $0^{\circ}$          & $30^{\circ}$     & $12$             & $0^{\circ}$      & $30^{\circ}$  & $6$            & $\hat{\theta},\hat{\phi}$  & $4$   \\   
    2 & $0^{\circ}$          & $30^{\circ}$     & $12$             & $0^{\circ}$      & $30^{\circ}$  & $6$            & $\hat{\theta},\hat{\phi}$  & $5$   \\
    3 & $0^{\circ}$          & $30^{\circ}$     & $12$             & $0^{\circ}$      & $30^{\circ}$  & $6$            & $\hat{\theta},\hat{\phi}$  & $6$   \\  
    4 & $0^{\circ}$          & $10^{\circ}$     & $36$             & $0^{\circ}$      & $10^{\circ}$  & $18$          & $\hat{\theta},\hat{\phi}$  & $4$   \\
    5 & $0^{\circ}$          & $10^{\circ}$     & $36$             & $0^{\circ}$      & $10^{\circ}$  & $18$          & $\hat{\theta},\hat{\phi}$  & $5$   \\
    6 & $0^{\circ}$          & $10^{\circ}$     & $36$             & $0^{\circ}$      & $10^{\circ}$  & $18$          & $\hat{\theta},\hat{\phi}$  & $6$   \\   
     \hline
  \end{tabular}
  \label{table_condition_Cube_12_12_12_4mm_Array_8x2x2_15mm_space}
\end{table}

\begin{figure}[!t]
  \centering
    \includegraphics[keepaspectratio, scale=0.55]{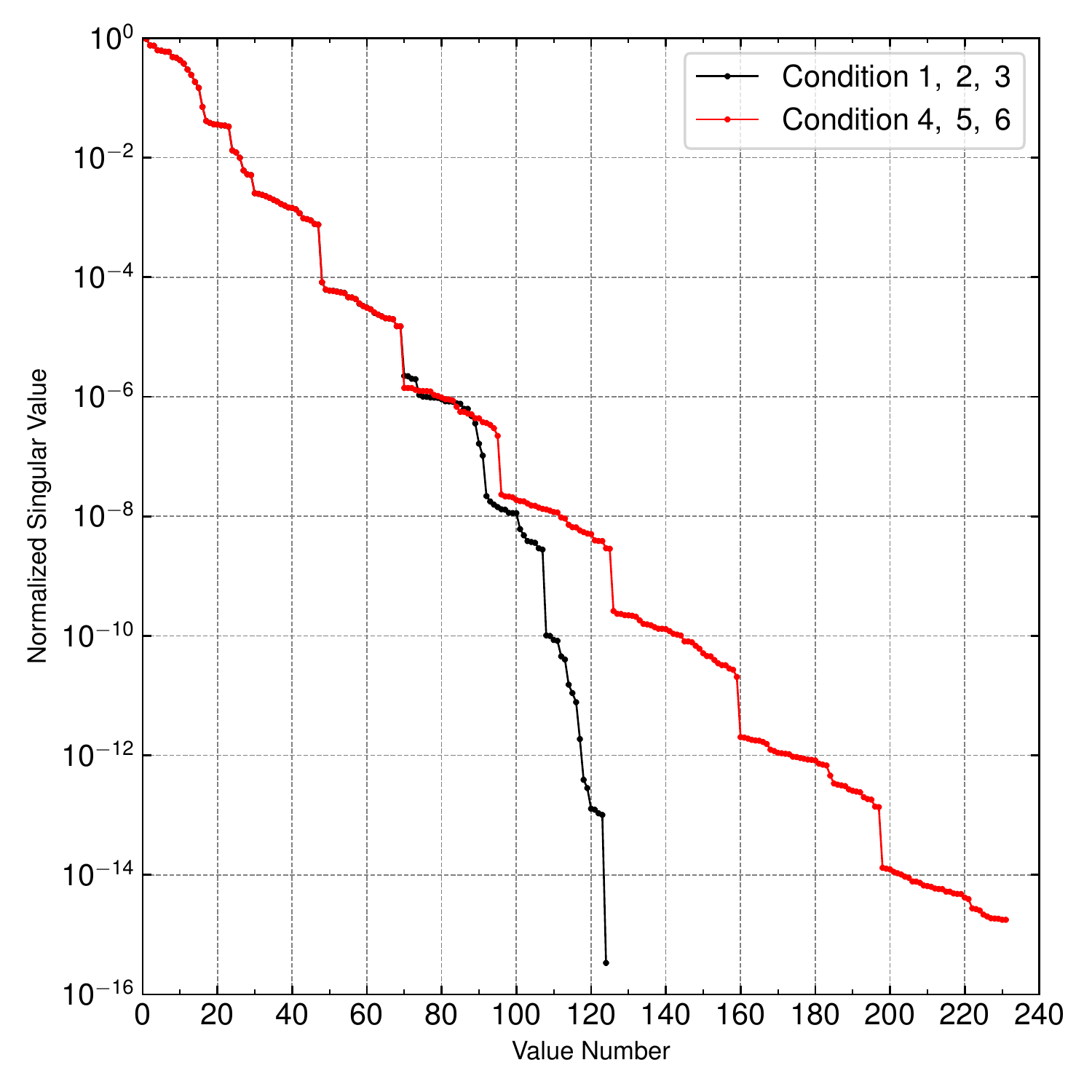}
  \caption{Normalized singular value of the cube.}
  \label{fig_svd_Cube_12_12_12_2mm_Array_8x2x2_15mm_space}
\end{figure}

\begin{figure}[!t]
  \begin{tabular}{c}
    \begin{minipage}[t]{0.45\hsize}
      \centering
      \includegraphics[keepaspectratio, scale=0.27]{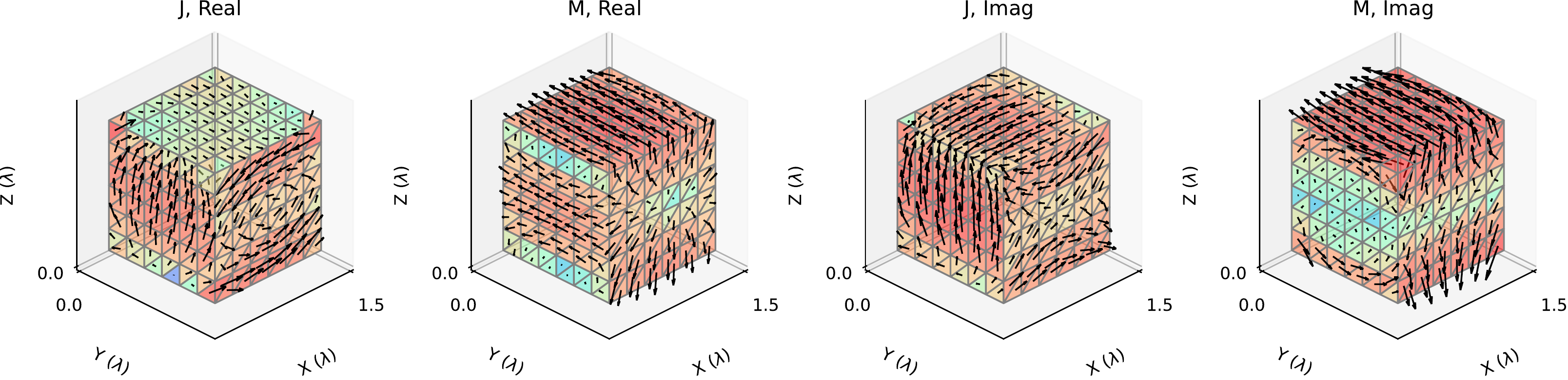}
      \subcaption{CBF for singular value \#$1$}
      \label{fig_Cube_12_12_12_2mm_Array_8x2x2_15mm_space_current0_CBF0_az45_el30}
    \end{minipage} \\
    \begin{minipage}[t]{0.45\hsize}
      \centering
      \includegraphics[keepaspectratio, scale=0.27]{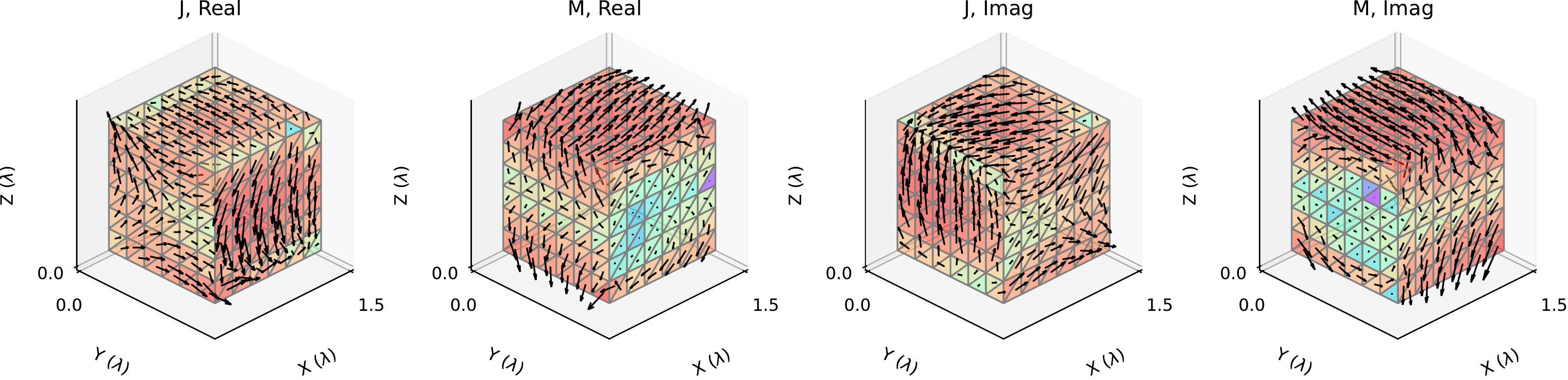}
      \subcaption{CBF for singular value \#$2$}
      \label{fig_Cube_12_12_12_2mm_Array_8x2x2_15mm_space_current0_CBF0_az45_el30}
    \end{minipage} \\
    \begin{minipage}[t]{0.45\hsize}
      \centering
      \includegraphics[keepaspectratio, scale=0.27]{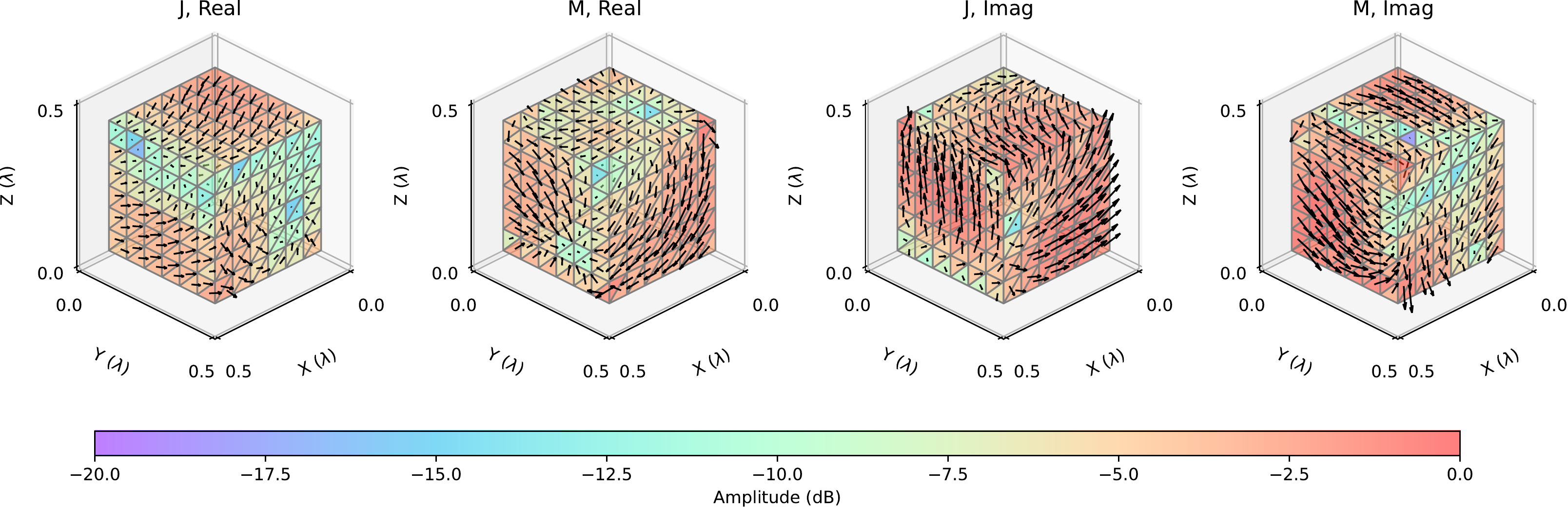}
      \subcaption{CBF for singular value \#$96$}
      \label{fig_Cube_12_12_12_2mm_Array_8x2x2_15mm_space_current0_CBF0_az45_el30}
    \end{minipage} 
  \end{tabular}
  \caption{CBF distributions for cell $1$ of the cube array. }
  \label{fig_Cube_12_12_12_2mm_Array_8x2x2_15mm_space_current0}
\end{figure}

\figuref\ref{fig_svd_Cube_12_12_12_2mm_Array_8x2x2_15mm_space} shows the distribution of singular values in \eqref{eq_gram_CBF_SVD} for a single cube. 
Note that these singular values are the same for all cells because the incident fields on each cell are identical except for the phase factor.
In this figure we can identify a few groups of singular values, although not as clearly as in the case of the sphere.
We see that the numbers of singular values in groups with small $l$s are the same as those in the corresponding sphere case (See \figuref\ref{fig_svd_Sphere_R5mm_Tess2}). 
Also, the singular values of groups with larger $l$ can be calculated more accurately with smaller incident angle intervals. 
\figuref\ref{fig_Cube_12_12_12_2mm_Array_8x2x2_15mm_space_current0} shows the real and imaginary parts of CBFs at the barycentric coordinate of each triangle in a cell represented by \eqref{eq_cbf_JM} for the cube array under the condition $6$ ($l=6$) in \tableref\ref{table_condition_Cube_12_12_12_4mm_Array_8x2x2_15mm_space}. 
The orthogonality between electric and magnetic CBFs is visible even with the CBF No.~$96$, which is the maximum number of CBFs to be considered with $l=6$. 

\begin{figure}[!t]
  \centering
    \includegraphics[keepaspectratio, scale=0.55]{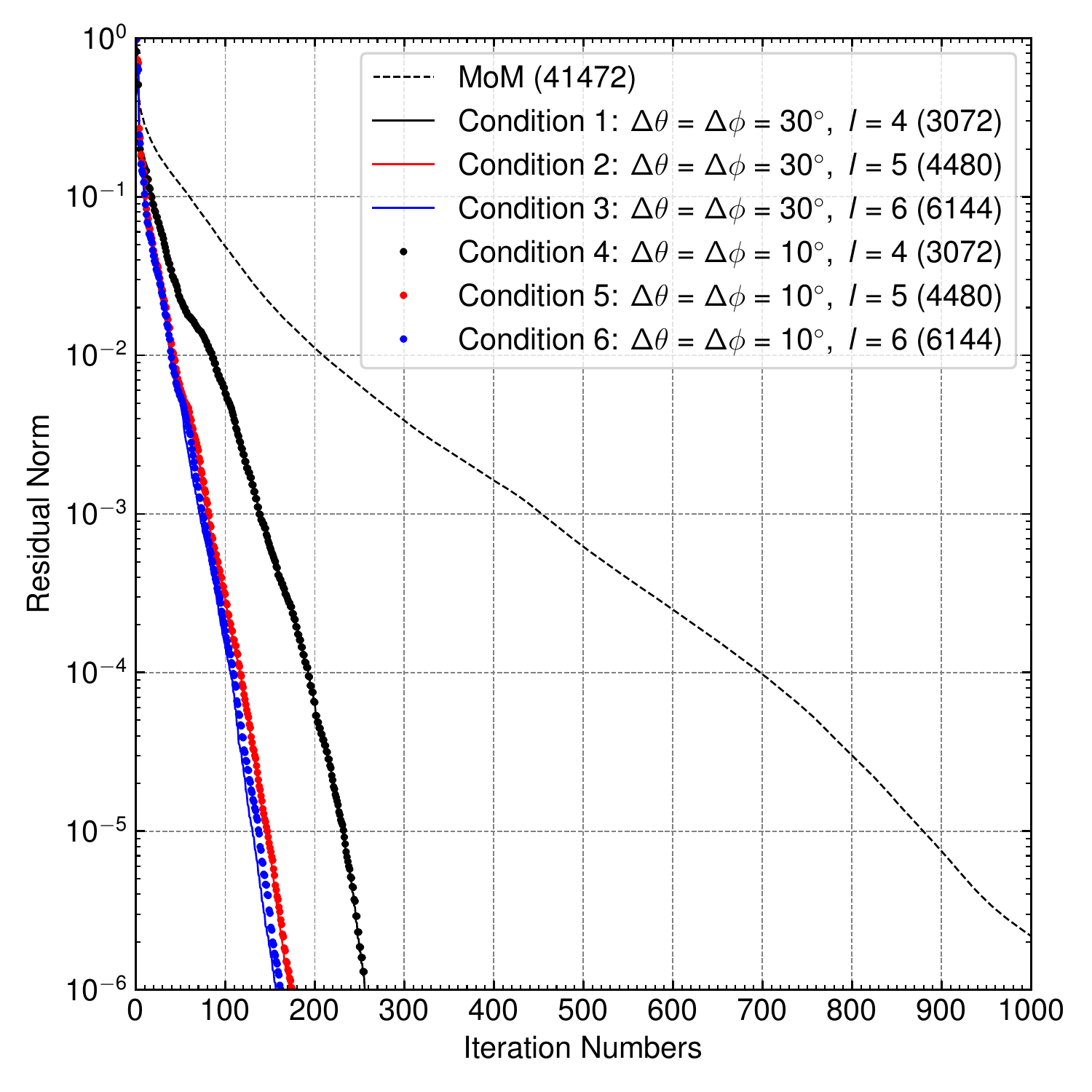}
  \caption{Convergence of the cube array analysis. Number in brackets of each series in the figure represents the unknowns of MoM and CBFM. }
  \label{fig_norm_Cube_12_12_12_2mm_Array_8x2x2_15mm_space}
\end{figure}

The relationship between the choice of CBF parameters, accuracy, and convergence can be interpreted as follows. 
First, we consider the relationship between the angular interval of the incident waves and accuracy. 
In this problem, the electric length corresponding to the circumference of the circumscribing sphere for a cube is $2\pi r \sqrt{\epsilon_{\rm r}} \simeq 3.8 \lambda$, where  $r = \sqrt{3} \lambda/5$.
Therefore, approximately $3.8$ cycles of electromagnetic currents could be generated on the surface of the fictitious circumscribing sphere, which tells that one may need at least $l\approx 4$ for an accurate resolution of the solution. 
On the other hand, the spacing of the incident wave may roughly correspond to the sampling interval of the electromagnetic current generated on the surface. 
This indicates that only $3$ points are sampled per cycle for the $l=4$ variation when the incident field interval is set to $30^{\circ}$. 
When the incident field interval is set to $10^{\circ}$, however, the number of sampling points per cycle for $l=4$ is raised to $9$ which is more reasonable. 
Next, we consider the relationship between the convergence of CBFM and the appropriate group number $l$ to be considered for CBFs. 
It is expected that the current distribution on a cube will be more complicated than that on a sphere due to corners and edges. 
This indicates that more CBFs are needed to represent complex currents on a cube than those determined by $l=4$.
From this consideration, we conclude that taking CBFs corresponding to $l$ greater than $5$ or $6$  will be desired in this problem.
The validity of this consideration is demonstrated in \figuref\ref{fig_norm_Cube_12_12_12_2mm_Array_8x2x2_15mm_space} which shows the convergence of the outer GMRES when the incident wave propagates in $-z$ direction.
The numbers in the brackets in the legend represent the number of unknowns of the MoM and the CBFM. 
In MoM, the residual norm $\delta_{\rm R}$ did not reach $10^{-6}$ even after $1000$ iterations. 
In CBFMs, on the other hand, this threshold is reached in less than $300$ iterations. 
We see that the convergence of CBFM is further improved when $l$ is set to $5$ or $6$ (condition $2$, $3$, $5$ and $6$), which is consistent with our consideration on the choice of CBFs based on the distribution of singular values in \eqref{eq_gram_CBF_SVD}.

\begin{figure}[!t]
  \centering
    \includegraphics[keepaspectratio, scale=0.55]{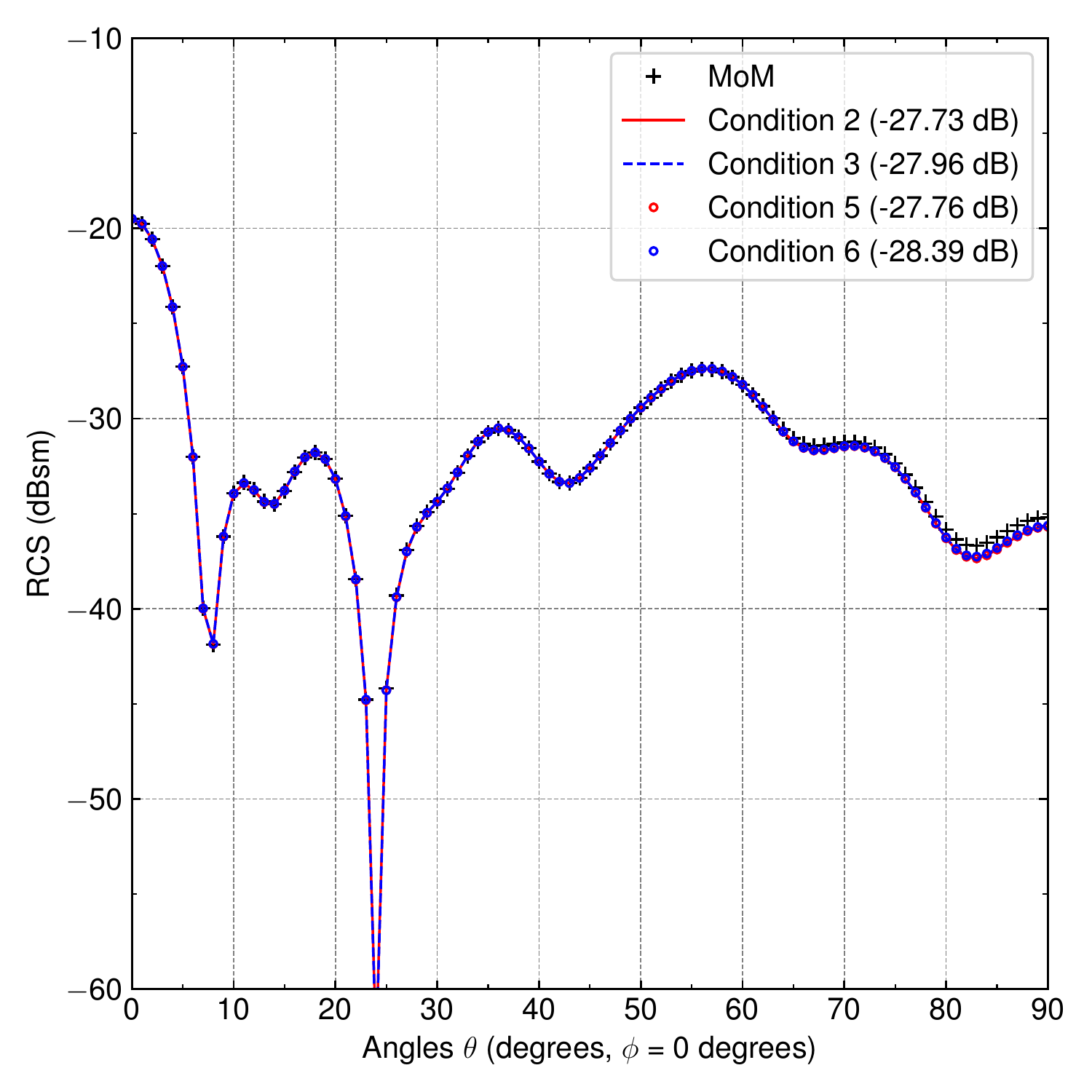}
  \caption{RCS of the cube array ($8\times2\times2$). The number in brackets of each series in the figure represents RMSE. }
  \label{fig_rcs_Cube_12_12_12_4mm_Array_8x2x2_15mm_space}
\end{figure}

\figuref\ref{fig_rcs_Cube_12_12_12_4mm_Array_8x2x2_15mm_space} shows the RCS patterns obtained with MoM, and the proposed method with conditions $2$, $3$, $5$, and $6$, where the criterion of the convergence is set as $\delta_{\rm R}<10^{-4}$.
The angle range, the angle resolution, and the polarization of the RCS pattern are the same as those for the sphere array analysis in \figuref\ref{fig_plot_rcs_Sphere_R5mm_Tess2_Array_8x16x2_15mm_space}.
The numbers in the brackets in the legend in \figuref\ref{fig_rcs_Cube_12_12_12_4mm_Array_8x2x2_15mm_space} represent the RMSE. 
We see that setting small values to angular intervals $\Delta\theta$ and $\Delta\phi$ slightly improves the accuracy of the analysis as mentioned in the previous discussion.
All these results show that the proposed method enables one to obtain accurate solutions with a small number of iterations as one considers singular values belonging up to $l=5$ or $6$ groups. 

The computational times for generating the primary CBFs and the final iteration for the condition $6$, relative to the total time taken to compute RCS pattern for $91$ directions ($0^{\circ}\leq\theta\leq90^{\circ}$, $\Delta\theta=1^{\circ}$) with MoM,
are less than $0.01$ and $0.21$, respectively. 
We thus see that the proposed method can compute the solution approximately five times faster than the conventional method in this example. 

These results show that the proposed method is applicable not only to smooth scatterers but also to scatterers with corners and edges, and can obtain accurate solutions faster than MoM. 
\subsection{Cylinder} \label{sec_Cylinder}
\begin{figure}[!t]
    \includegraphics[keepaspectratio, scale=0.35]{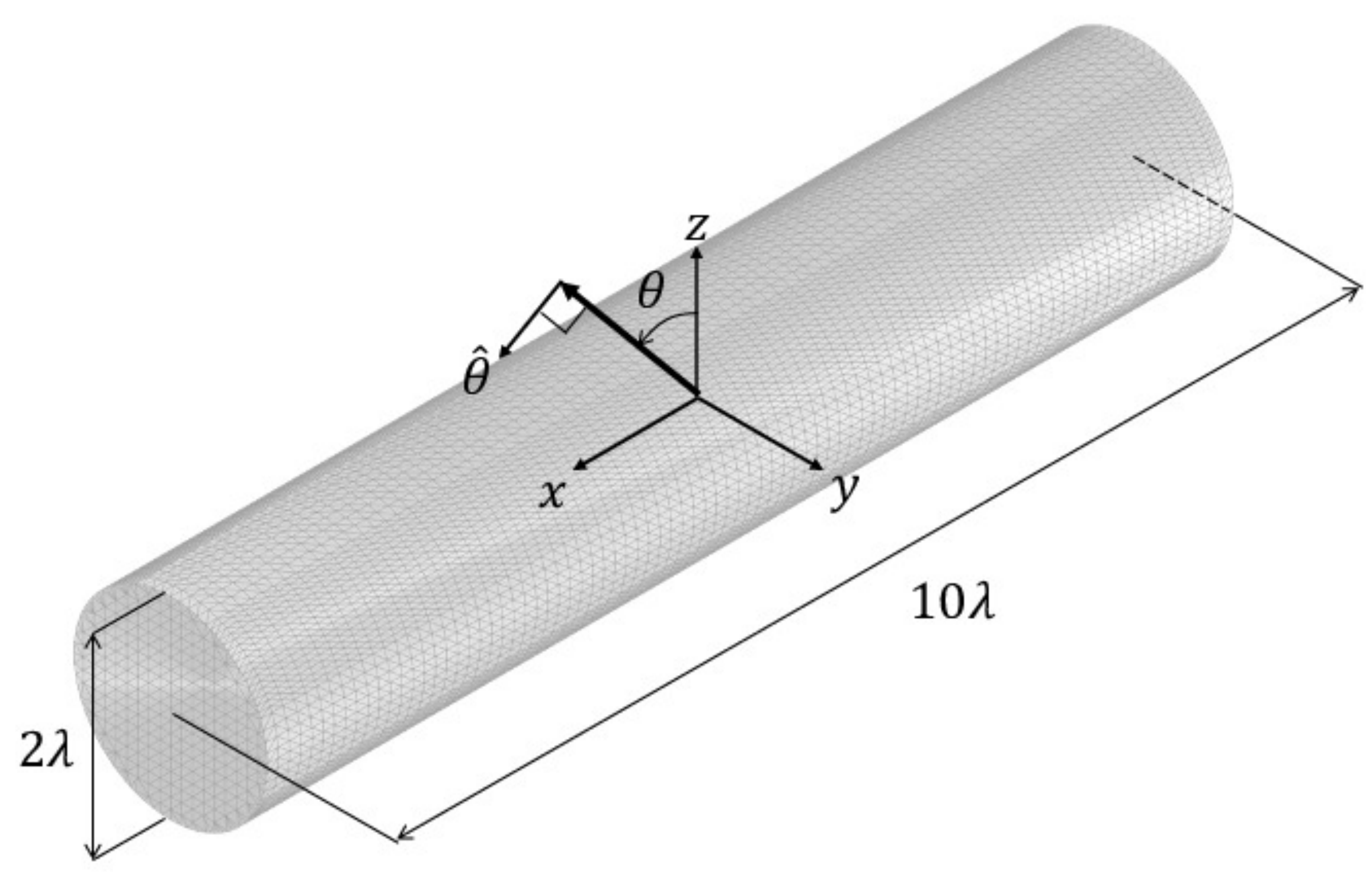}
  \caption{Cylinder}
  \label{fig_cylinder}
\end{figure}

Finally we analyze the cylinder shown in \figuref\ref{fig_cylinder} as an example of a large connected scatterer. 
The angle range and polarization to calculate the RCS pattern are $0^{\circ} \leq \theta \leq 45^{\circ}$, $\Delta\theta=0.5^{\circ}$, $\phi=0^{\circ}$ and $\hat{\theta}$ ($91$ directions). 
The number of unknowns in MoM is $52752$. 

We use IPCBFs as the CBFs in this analysis since this scatterer cannot be split into more than one closed scatterers by any cell arrangement.
\tableref\ref{table_condition_cylinder2} shows the conditions for the IPCBF generation. 
The plane waves used for the IPCBF generation are with $\hat{\theta}$ polarization and their propagation directions are taken in the same coordinate plane ($\phi=0^{\circ}$ plane) on which the RCS is calculated. 
The number of plane wave directions for IPCBF generation $s$ is equal to $N_{\theta} (=10)$ in all the conditions.
The allowance of the residual norm $\delta_{\rm r}$ (see \ref{sec_CBF_Generation}) used to determine IPCBFs is set to $1.0\times10^{-1}$ to $1.0\times10^{-3}$. 
The cells defining CBFs are cubes with the side lengths of $1.25\lambda$, and the number of cells is $M=32$. 
The threshold $\delta_{\rm SVD}$ for SVD (see the paragraph below \eqref{eq_gram_CBF_SVD}) is set as $1.0\times10^{-3}$ and $1.0\times10^{-7}$. 
For the purpose of comparison, we have also carried out the same analysis using the primary CBFs generated using parameters in
\tableref\ref{table_condition_cylinder_pcbf}.
We use CMP in all these analyses.
\begin{table}[!t]
  \caption{Parameters for Calculating IPCBFs for Cylinder}
  \centering
  \begin{tabular}{cccccccccc} \hline
    No. & $\delta_{\rm r}$    & $\theta_{\rm s}$ & $\Delta\theta$ & $N_{\theta}$ & Cell ($\lambda$) & $\delta_{\rm SVD}$ \\ \hline \hline
    $1$ & $1.0\times10^{-3}$  & $0^{\circ}$      & $5^{\circ}$    & $10$         & $5/4$                & $1.0\times10^{-3}$ \\
    $2$ & $1.0\times10^{-2}$  & $0^{\circ}$      & $5^{\circ}$    & $10$         & $5/4$                & $1.0\times10^{-3}$ \\ 
    $3$ & $1.0\times10^{-1}$  & $0^{\circ}$      & $5^{\circ}$    & $10$         & $5/4$                & $1.0\times10^{-3}$ \\ 
    $4$ & $1.0\times10^{-3}$  & $0^{\circ}$      & $5^{\circ}$    & $10$         & $5/4$                & $1.0\times10^{-7}$ \\ \hline
  \end{tabular}
  \label{table_condition_cylinder2}
\end{table}
\setlength{\tabcolsep}{1.5mm}
\begin{table}[!t]
  \caption{Parameters for Calculating primary CBFs for Cylinder}
  \centering
  \begin{tabular}{ccccccccc} \hline
    $\theta_{\rm s}$  & $\Delta\theta$ & $N_{\theta}$ & $\phi_{\rm s}$ & $\Delta\phi$ & $N_{\phi}$ & Pol.                                & Cell ($\lambda$) & $\delta_{\rm SVD}$\\ \hline \hline
    $0^{\circ}$          & $10^{\circ}$     & $36$               & $0^{\circ}$      & $10^{\circ}$  & $18$            & $\hat{\theta},\hat{\phi}$  & $5/4$               & $1.0\times10^{-3}$ \\ \hline
  \end{tabular}
  \label{table_condition_cylinder_pcbf}
\end{table}

\begin{figure}[!t]
    \includegraphics[keepaspectratio, scale=0.55]{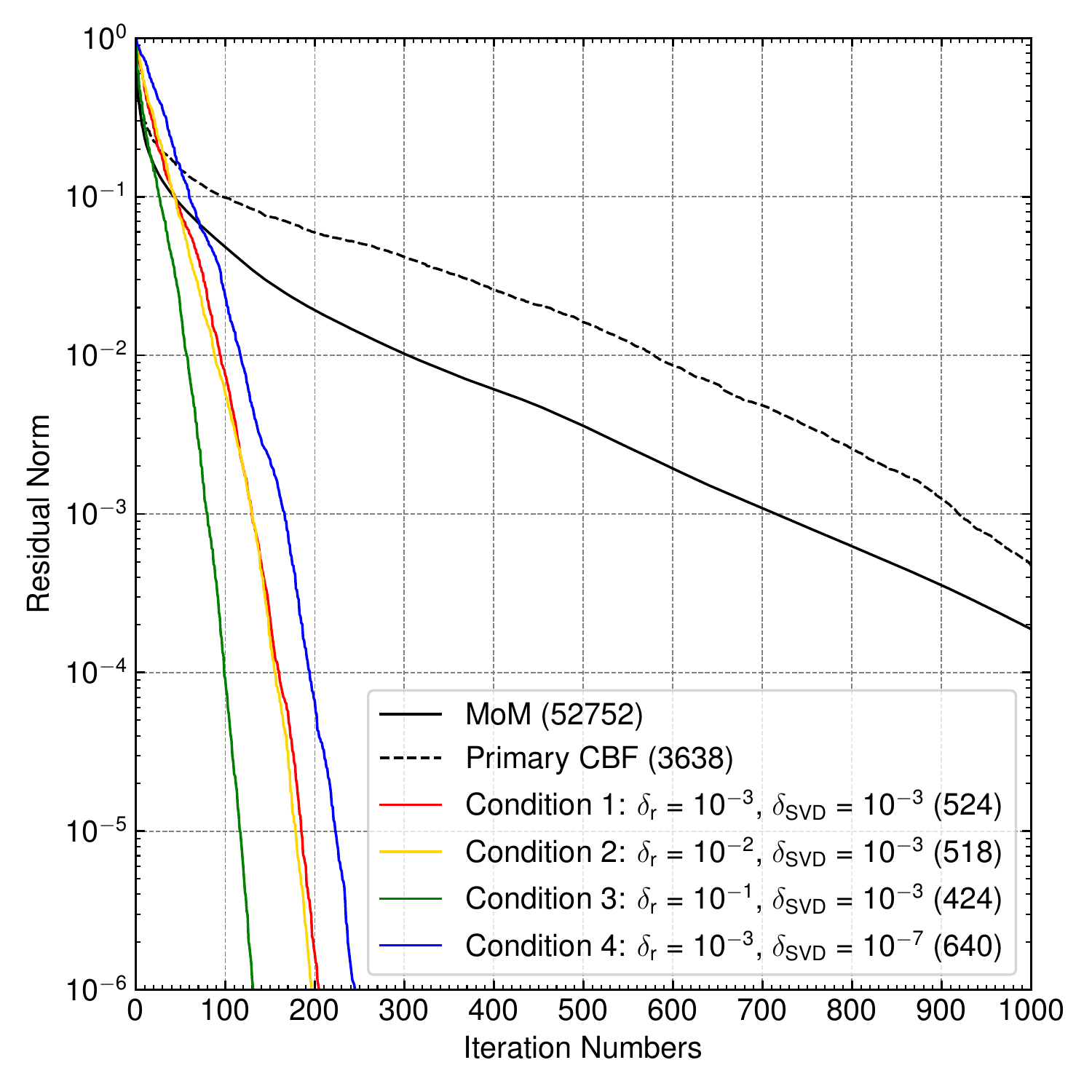}
  \caption{Convergence of the cylinder analysis for the incident wave with $\hat{\theta}$--polarization from $-z$ direction. Number in brackets of each series in the figure represents the unknowns of MoM and CBFM.}
  \label{fig_norm_cylinder}
\end{figure}

\figuref\ref{fig_norm_cylinder} shows the convergence property of the outer GMRES when the $\hat{\theta}$ polarized incident wave propagates in $-z$ direction. 
The numbers in the brackets in the legend represent the number of unknowns of the MoM and the CBFM. 
The CBFM with IPCBFs converges much faster than MoM.
In contrast to the previous analyses, however, the CBFM with the primary CBFs shows even poorer convergence compared to that of the MoM. 
This is thought to be due to the effect of the fictitious open edges that occur when the scatterer is divided in the CBF generation process. 
As a matter of fact, it is not reasonable to solve problems for a dielectric object containing open edges with the PMCHWT formulation, as noted in section~\ref{sec_CBF_Generation}. 
This problem does not arise when the scatterer can be divided into disjoint closed scatterers, as in examples in the previous sections. 

\begin{figure}[!t]
  \centering
    \includegraphics[keepaspectratio, scale=0.55]{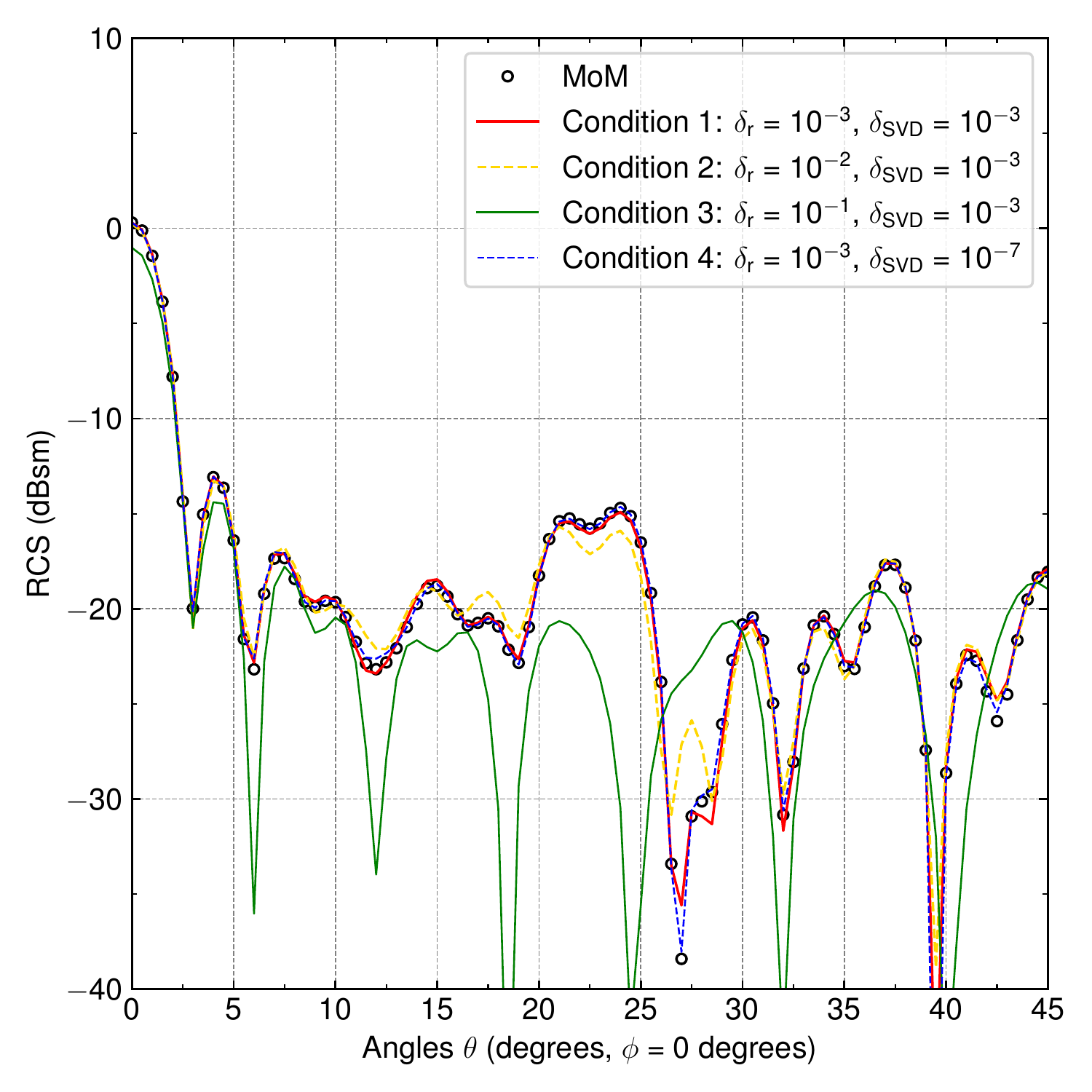}
  \caption{RCS of the cylinder}
  \label{fig_rcs_cylinder2_e-04}
\end{figure}

\begin{figure}[!t]
  \centering
    \includegraphics[keepaspectratio, scale=0.6]{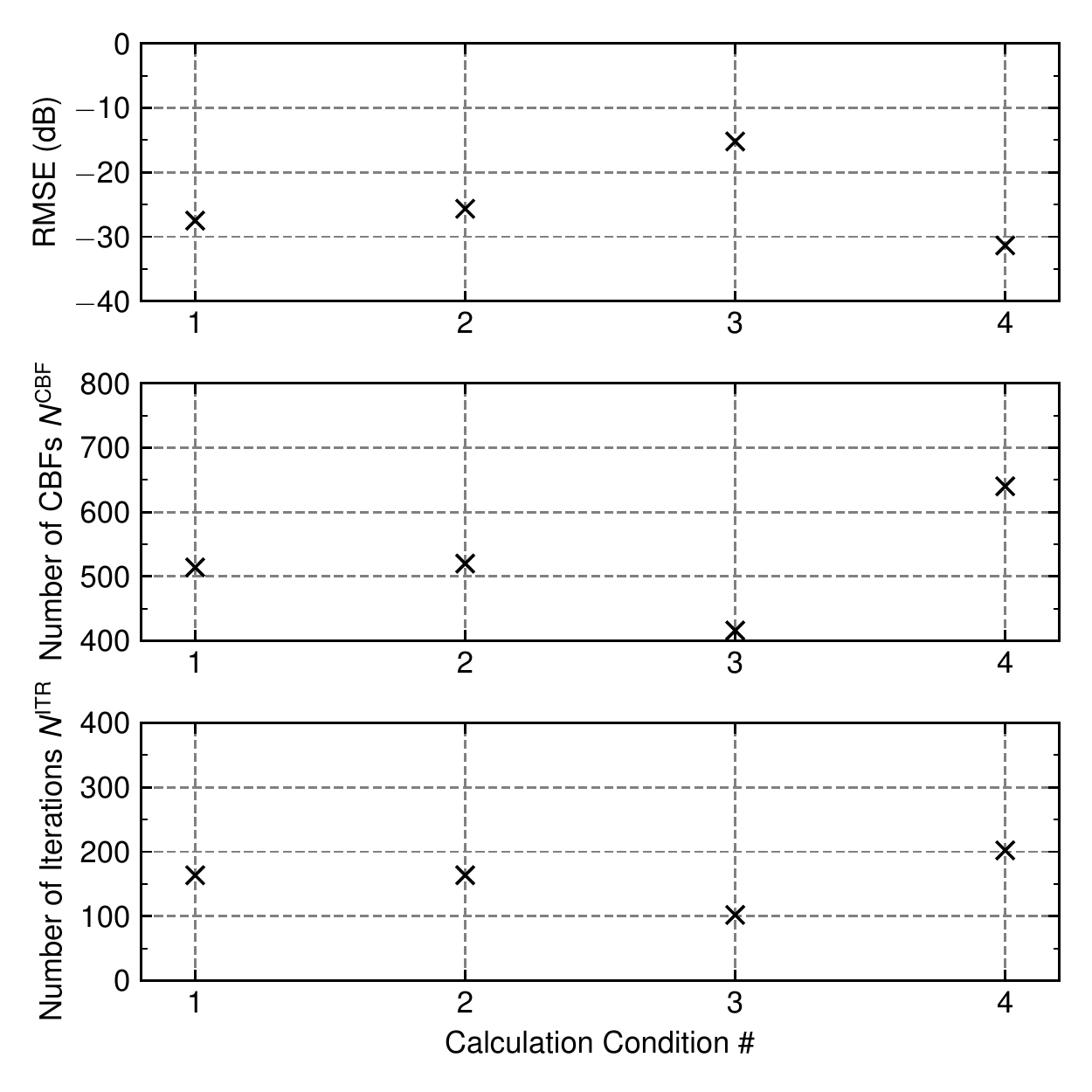}
  \caption{Accuracy, number of CBFs, and average number of final iterations}
  \label{fig_rmse_cylinder}
\end{figure}

\figuref\ref{fig_rcs_cylinder2_e-04} and \figuref\ref{fig_rmse_cylinder} show the RCS pattern, accuracy (RMSE), number of CBFs $N^{\rm CBF}$, and average number of final iterations $N^{\rm ITR}$ of the CBFM which uses IPCBFs. 
The number $N^{\rm ITR}$ gives the mean value of the numbers of iterations in the analyses for $91$ incident directions. 
The number $N^{\rm CBF}$ does not change significantly with conditions and the average numbers $N^{\rm ITR}$ for conditions 1 and 2 are almost the same.
The RCS obtained with the third condition does not agree with the MoM result, though this condition requires less $N^{\rm ITR}$ than the first and second conditions.
The RCS obtained with the second condition is better than the result with the third condition and is reasonably accurate (RMSE: $-24.26$dB). 
The RCS obtained with the first condition agrees well with the MoM results down to low levels and RMSE is $-26.49$dB. 

\tableref\ref{table_computational_time_cylinder2} shows the relative computational time ratio, where
``CBF Gen.", ``Iter." and ``Total" represent the relative computational time for CBF generation, iteration for RCS analysis, and their total relative to the total computational time to analyze cases for $91$ incident directions with the MoM, respectively.
The generation time of CBFs for the second condition is shorter than that for the first one. 
In both cases, the solution is obtained more than five times faster than in the conventional MoM. 
These results indicate that there is a trade-off between accuracy and analysis time. 
To obtain a more accurate solution, the value of the residual norm $\delta_{\rm r}$ should be around $1.0 \times 10^{-3}$, while
this value should be around $1.0 \times 10^{-2}$ for a rough and fast study of the overall trend. 

In the fourth condition, the threshold of singular value decomposition $\delta_{\rm SVD}$ (see \ref{sec_orthogonalization}) is reduced while other parameters are kept the same as in the first condition.
We see that the accuracy further improves so that the RMSE becomes $-31.35$dB. 
This result tells that we can further control the accuracy of the analysis by reducing $\delta_{\rm SVD}$ when the the residual norm $\delta_{\rm r}$ is less than $1.0 \times 10^{-3}$.
\begin{table}[!t]
  \caption{Relative Computational Time Ratio of Cylinder}
  \centering
  \begin{tabular}{ccc|c} \hline
    No. & CBF Gen. & Iter.  & Total   \\ \hline \hline
    $1$ & $0.06$   & $0.12$ & $0.18$  \\
    $2$ & $0.02$   & $0.12$ & $0.14$  \\ 
    $3$ & $<0.01$  & $0.08$ & $<0.09$ \\ 
    $4$ & $0.06$   & $0.16$ & $0.22$  \\ \hline
  \end{tabular}
  \label{table_computational_time_cylinder2}
\end{table}
\section{Conclusion}\label{sec_conclusion}
In this paper, we proposed a CBFM and a CMP which use the duality of electromagnetic currents.  These techniques provide an efficient method to analyze scattering problems for homogeneous dielectric materials. 
The numerical results for various scatterers confirm the effectiveness of the proposed method. 
Furthermore, we investigated mathematically the relationship between the distribution of singular values of a certain Gram matrix and the spacing of incident angles for plane waves for the generation of CBFs for disjoint scatterers.
We also carried out parametric studies for the relationship between the accuracy of the solution and the threshold of singular values of a certain Gram matrix  for connected scatterers.

We may mention the following as possible future works of this research:
The proposed CMP is not very efficient for a connected scatterer while it can significantly accelerate convergence in multiple scatter problems as mentioned in section \ref{sec_Cylinder}. 
The use of the  proposed method basically improves the convergence even in connected scatterers, as has been described in this paper. 
However, one may make the analyses in this case faster by orthogonalizing CBFs for electric and magnetic components taking into account the effect of the fictitious open edges in (\ref{eq_gram_CBF})--(\ref{eq_CBF}). The details of such approach remain to be investigated.
The second prospect is about the extension of the proposed method to more complex media. 
The proposed method has the potential to be extended to the analysis of more complex structures, such as scatterers containing multiple objects with different dielectric constants and conductors. 
The extension to the former problem would be relatively straightforward, but the case containing conductors would be more complicated. 
We are currently working on these problems. 


%


\ifCLASSOPTIONcaptionsoff
  \newpage
\fi


\begin{thebibliography}{99}
\bibitem{Harrington1993}
  R.~F.~Harrington, \emph{Field computation by the moment methods}, Wiley-IEEE Press, New York, 1993. 
\bibitem{Saad2003}
  Y.~Saad, \emph{Iterative methods for sparse linear systems}, SIAM, 2003. 
\bibitem{Chew2001}
  W.~C.~Chew, J.~M.~Jin, E.~Michielssen, J.~Song \emph{Fast and efficient algorithms in computational electromagnetics}, Altech House, Boston, 2001. 
\bibitem{Ergul2014}
  \"{O}~Erg\"{u}l and L.~G\"{u}rel, \emph{The multilevel fast Multipole algorithm (MLFMA) for solving large-scale computational electromagnetics problems}, Wiley-IEEE Press, New York, 2014. 
\bibitem{Song1997}
    J.~Song, C.~-C.~Lu, and W.~C.~Chew, ``Multilevel fast multipole algorithm for electromagnetic scattering by large complex objects,'' \emph{IEEE Trans. Antennas Propag.}, vol.~45, no.~10, pp.~1488--1493, Oct. 1997. 
\bibitem{Grasedyck2003}
L.~Grasedyck, and W.~Hackbusch, ``Construction and arithmetics of $\mathcal{H}$-matrices,'' \emph{Computing}, vol.~70, pp.~295--334, Jul. 2003.
\bibitem{Heldrig2007}
  A.~Heldrig, J.~M.~Rius, J.~M.~Tamayo, J.~Parr\'{o}n, E.~\'{U}beda, ``Fast direct solution of method of moments linear system,'' \emph{IEEE Trans. Antennas Propag.}, vol.~55, no.~11, pp.~3220--3228, Nov. 2007. 
\bibitem{Gholami2018}
  R.~Gholami, J.~Mojolagbe, A.~Menshov, F.~S.~H.~Lori, and V.~Okhmatovski, ``$\mathcal{H}$-matrix arithmetic for fast direct and iterative method of moment solution of surface-volume-surface EFIE for 3-D radiation problems,'' \emph{Progress In Electromagnetics Research B}, vol.~82, pp.~189--210, Dec. 2018. 
\bibitem{Andriulli2008}
  F.~P.~Andriulli, K.~Cools, H.~Ba\u{g}ci, F.~Olyslager, A.~Buffa, S.~Christiansen and E.~Michielssen, ''A multiplicative Calder\'{o}n preconditioner for the electric field integral equation,'' \emph{IEEE Trans. Antennas Propag.}, vol.~56, no.~8, pp.~2398--2412, Aug. 2008.
\bibitem{Rao1982}
  S.~Rao, D.~Wilton, and A.~Glisson, ``Electromagnetics scattering by surfaces of arbitrary shape,'' \emph{IEEE Trans. Antennas Propag.}, vol.~30, no.~3, pp.~409--418, May 1982.
\bibitem{Buffa2007}
  A.~Buffa and S.~H.~Christiansen ``A dual finite element complex on the barycentric refinement,'' \emph{Mathematics of Computation}, vol.~76, no.~260, pp.~1743--1769, Oct. 2007.
\bibitem{Yan2010}
  S.~Yan and J.~M.~Jin and Z.~Nie, ``A comparative study of Calder\'{o}n preconditioners for PMCHWT equations,'' \emph{IEEE Trans. Antennas Propag.}, vol.~58, no.~7, pp.~2375--2383, Oct. 2010.
\bibitem{Cools2011}
  K.~Cools, F.~P.~Andriulli and E.~Michielssen, ``A Calder\'{o}n multiplicative preconditioner for the PMCHWT integral equation,'' \emph{IEEE Trans. Antennas Propag.}, vol.~59, no.~12, pp.~579--4587, Dec. 2011.
\bibitem{Niino2012}
  K.~Niino and N.~Nishimura, ``Calder\'{o}n preconditioning approaches for PMCHWT formulations for Maxwell's equations,'' \emph{International Journal of Numerical Modelling: Electronic Networks, Devices and Fields}, vol.~25, no.~5--6, pp.~558--572, Sep. 2012.
\bibitem{yla2012stable}
  P.~Yla-Oijala, S.~P.~Kiminki, K.~Cools, F.~P.~Andriulli and S.~Jarvenpaa, ``Stable discretization of combined source integral equation for scattering by dielectric objects,'' \emph{IEEE Trans. Antennas Propag.}, vol.~60, no.~5, pp.~2575--2578, May 2012.
\bibitem{Matekovits2007}
  L.~Matekovits, V.~A.~Laza, and G.~Vecchi, ``Analysis of large complex structures with the synthetic-functions Approach,'' \emph{IEEE Trans. Antennas Propag.}, vol.~55, no.~9, pp.~2509--2521, Sep. 2007.
\bibitem{Prakash2003}
  V.~Prakash, and R.~Mittra, ``Characteristic basis function method: a new technique for efficient solution of method of moments matrix equations,'' \emph{Microw. Opt. Techn. Let.}, vol.~36, no.~2, pp.~95--100, Jan. 2003. 
\bibitem{Lucente2008}
  E.~Lucente, A.~Monorchio, and R.~Mittra, ``An iteration-free MoM approach based on excitation independent characteristic basis functions, for solving large multiscale electromagnetic scattering problems,'' \emph{IEEE Antennas Propag.}, vol.~56, no.~4, pp.~999--1007, Apr. 2008. 
\bibitem{Mittra2008}
  R.~Mittra and K.~Du, ``Characteristic basis function method for iteration-free solution of large method of moments problems,'' \emph{Progress In Electromagnetics Research B}, vol.~6, pp.~307--336, Apr. 2008. 
\bibitem{Maaskant2008}
  R.~Maaskant, R.~Mittra, and A.~G.Tijhuis. ``Fast analysis of large antennas arrays using the characteristic basis function method and the adaptive cross approximation algorithm,'' \emph{IEEE Trans. Antennas Propag.}, vol.~56, no.~11, pp.~3440--3451, Nov. 2008.
\bibitem{Garcia2008}
  E.~Garcia, C.~Delgado, I.~G.~Diego, and M.~F.~Catedra, ``An iterative solution for electrically large problems combining the characteristic basis function method and the multilevel fast multipole algorithm,'' \emph{IEEE Trans. Antennas Propag.}, vol.~56, no.~8, pp.~2363--2371, Aug. 2008.
\bibitem{Chen2015}
X.~Chen, C.~Gu, J.~Ding, Z.~Li, and Z.~Niu, ``Multilevel fast adaptive cross-approximation algorithm with characteristic basis functions,'' \emph{IEEE Trans. Antennas Propag.}, vol.~63, no.~9, pp.~3994--4002, Sep. 2015. 
\bibitem{Tanaka2016}
  T.~Tanaka, Y.~Inasawa, Y.~Nishioka, H.~Miyashita, ``Improved primary characteristic basis function method for monostatic radar cross section of specific coordinate plane,'' \emph{IEICE Trans. Electron.}, vol.~99-C, no.~1, pp.~28--35, Jan. 2016.
\bibitem{Tanaka2017}
  T.~Tanaka, Y.~Inasawa, Y.~Nishioka, H.~Miyashita, ``Improved primary characteristic basis function method considering higher-order multiple scattering,'' \emph{IEICE Trans. Electron.}, vol.~100-C, no.~1, pp.~45--51, Jan. 2017.
\bibitem{Tanaka2019}
  T.~Tanaka, K.~Niino, N.~Nishimura, M.~Takikawa, and N.~Yoneda, ``A generation scheme of the characteristic basis functions by using block Krylov subspace algorithm (Japanese),'' \emph{Trans. Jpn. Soc. Comput. Methods Eng.}, vol.~19, pp.~99--102, Dec. 2019. 
\bibitem{Knott1993}
  E.~F.~Knott, J.~F.~Shaeffer, and M.~T.~Tuley, \emph{Radar cross section second edition}, Artech House, Boston, 1993.
\bibitem{niino-jcp}
  K.~Niino, and N.~Nishimura,  ``Preconditioning based on Calderon's formulae for periodic fast multipole methods for Helmholtz' equation,'' \emph{J. Comp. Phys.}, vol.~231, pp.~66--81, Jan. 2012.
\end{thebibliography}
\end{document}